\documentclass[11pt,psamsfonts]{amsart}
\usepackage{amssymb,amscd,amsmath,amsthm,mathscinet}

\input prepictex
\input pictex
\input postpictex

\chardef\bslash=`\\ 

 \theoremstyle{plain} 
\newtheorem{theorem}{Theorem}[section]
\newtheorem{corollary}[theorem]{Corollary}
\newtheorem{lemma}[theorem]{Lemma}
\newtheorem{proposition}[theorem]{Proposition}

\theoremstyle{remark}
\newtheorem{remark}[theorem]{Remark}
\newtheorem{example}[theorem]{Example}
\theoremstyle{definition}
\newtheorem{definition}[theorem]{Definition}

\newcommand{\thmref}[1]{Theorem~\ref{#1}}

\newcommand{\proref}[1]{Proposition~\ref{#1}}
\newcommand{\lemref}[1]{Lemma~\ref{#1}}

\def\idx{\operatorname{id}}
\def\coker{\operatorname{coker}}
\def\ker{\operatorname{ker}}
\def\N{{\mathbb N}}

\def\Q{{\mathbb Q}}
\def\C{{\mathbb C}}

\def\NN{{{\mathbb N}}}

\def\Z{{\mathbb Z}}

\def\udodot{\cdot}

\newcommand{\End}{\operatorname{End}}
\newcommand{\Aut}{\operatorname{Aut}}

\def\inv{^{-1}}

\numberwithin{equation}{section}

\def\Splited#1>#2#3<#4{\unitlength.65mm%
\mbox{$\displaystyle #1$\rule[-3.6mm]{0mm}{3.6mm}%
\begin{picture}(25,10)(0,5.5)
      \put(2.5,7){\vector(1,0){20}}
      \put(12,8.5){$\scriptstyle #2$}
      \put(17.5,5){\oval(10,10)[br]}
      \put(17.5,0){\vector(-3,1){15}}
      \put(12.5,2){$\scriptstyle #3$}
\end{picture}
$\displaystyle #4$}}
\def\picArrow#1{\unitlength.65mm%
\begin{picture}(25,10)(0,5.5)
      \put(2.5,7){\vector(1,0){20}}
      \put(12,8.5){$\scriptstyle #1$}
\end{picture}}
\def\SES#1#2>#3>#4#5<#6{%
\mbox{$\displaystyle 1$\picArrow{}
$\displaystyle #1$\picArrow#2\Splited#3>#4#5<{#6}
\picArrow{}$\displaystyle 1$}}
\def\picarrow#1{\unitlength.35mm%
\begin{picture}(25,10)(0,9.3)
      \put(2.5,13){\vector(1,0){20}}
      \put(12,15){$\scriptstyle #1$}
\end{picture}}
\def\ses#1#2>#3>#4#5<#6{%
\mbox{$\textstyle 1$\picarrow{}
$\textstyle #1$\picarrow#2\splited#3>#4#5<{#6}
\picarrow{}$\textstyle 1$}}
\def\splited#1>#2#3<#4{\unitlength.35mm%
\mbox{$\textstyle #1$\rule[-3.00mm]{0mm}{3.00mm}%
\begin{picture}(25,10)(0,9.3)
      \put(2.5,13){\vector(1,0){20}}
      \put(12,15){$\scriptstyle #2$}
      \put(17.5,11){\oval(10,10)[br]}
      \put(17.5,6){\vector(-3,1){15}}
      \put(17.5,0){$\scriptstyle #3$}
\end{picture}
$\textstyle #4$}}

\begin{document}

\title{Hecke algebras from groups acting on trees \\ and HNN extensions}
\author[U. Baumgartner]{Udo Baumgartner}
\address{School of Mathematical and Physical Sciences,
              The University of Newcastle, University Drive, Building V,
              Callaghan, NSW 2308, Australia.}
\email{udo.baumgartner@newcastle.edu.au}
\author[M. Laca]{Marcelo Laca}
\address{Department of Mathematics, University of Victoria, Victoria, BC, Canada.}
\email{laca@math.uvic.ca}
\author[J. Ramagge]{Jacqui Ramagge}
\address{School of Mathematics and Applied Statistics, University of Wollongong, 
Wollongong NSW 2522, Australia.}
\email{ramagge@uow.edu.au}
\author[G. Willis]{George Willis}
\address{School of Mathematical and Physical Sciences,
              The University of Newcastle, University Drive, Building V,
              Callaghan, NSW 2308, Australia.}
\email{george.willis@newcastle.edu.au}
\thanks{Research supported by the Australian Research Council and the Natural Sciences and Engineering Research Council of Canada.}
\date\today

\begin{abstract} 
We study 
Hecke algebras of 
groups acting on 
trees 
with respect to 
geometrically defined 
subgroups. 
In particular, 
we consider 
Hecke algebras of 
groups of automorphisms of 
locally finite trees 
with respect to 
vertex 
and 
edge stabilizers 
and 
the stabilizer 
of an end 
relative to 
a vertex stabilizer, 
assuming that 
the actions 
are 
sufficiently transitive. 
We focus on 
identifying 
the structure of 
the resulting Hecke algebras, 
give explicit multiplication tables 
of the canonical generators 
and 
determine whether 
the Hecke algebra 
has  a universal $C$*-completion. 
The paper 
unifies 
past algebraic and analytic approaches 
by focusing on 
the common geometric thread.
The results 
have implications for 
the general theory 
of totally disconnected locally compact groups. 
\end{abstract}
\date{\today}
\maketitle

\section*{Introduction}
\label{sec:intro} 

In this paper 
we study 
Hecke algebras 
associated with
groups 
acting on 
locally finite trees. 
Hecke algebras 
are a representation-theoretic tool. 
Just as 
representations of 
a group $G$ 
on a complex vector space 
correspond to 
representations of 
its group algebra $\mathbb{C}[G]$, 
the algebra representations of 
the Hecke algebra $\mathbb{C}[G,G_0]$ 
of $G$ 
relative to 
an almost normal subgroup $G_0$ 
correspond to 
linear representations 
of $G$
which are 
generated by 
their $G_0$-fixed vectors. 
In the context of an arbitrary subgroup $G_0$, such representations are called 
{\em smooth} or \emph{algebraic}
representations of the group $G$, see
\cite{H-C.harm.anal,B.centre,BK.types} and~\cite[Chapter 1]{kri}. 
The groups $G$ considered here are subgroups of the automorphism group of a tree, with the almost normal subgroups $G_0$ being either the stabilizer of a vertex or the fixator of an edge of the tree.

Groups of tree automorphisms that do not fix an end of the tree are treated first. It is seen that in fact all groups $G$ of tree automorphisms that are sufficiently highly transitive yield the same Hecke algebra in each of the cases when $G_0$ is the stabilizer of a vertex and the fixator of an edge.
In the vertex stabilizer case, this Hecke algebra is the algebra of radial functions on the tree, as already calculated in~\cite{Macdonald_spherical} and~\cite{harmAna+repThry(Gs-act>homogTrees)}.

Groups of automorphisms
that do fix an end are of particular importance. 
These occur as 
HNN-extensions through Bass-Serre theory \cite{ser-tre} and, essentially for that reason, they play a role in the structure theory of totally disconnected locally compact groups analogous to that filled by the $(ax+b)$-group in Lie theory, see \cite[Proposition 2]{wil-str}, \cite{wil-comm} and the tree representation theorem \cite[Theorem 4.1]{contrG+scales(AUT(tdlcG))}. The representation theory of these groups is therefore basic to the representation theory of general totally disconnected groups. In this case the Hecke algebra depends very much on which subgroup $G$ of the full group of automorphisms fixing the end is considered. We initiate the study of these algebras by computing the structure of a few examples in the last section, 
giving explicit 
multiplication tables 
of the canonical generators. 

In the functional analytic setting, it is unitary representations on Hilbert space that are of importance and these correspond to representations of a $C$*-algebra. 
Inversion of group elements gives rise 
a *-algebraic structure on a Hecke algebra, and 
we are also able to decide in some of our examples whether 
the Hecke *-algebra
has an enveloping $C$*-algebra, a question that
has recently attracted 
a lot of interest, 
see%
~\cite{bin}, \cite{bos-con}, 
~\cite{hall, tza, primbc, lac-lar, etal, hecke-nonuniversal}. 
We use
general results on Hecke algebras 
for semidirect products to show that
our Hecke algebras 
and $C$*-algebras are semigroup crossed products. 
We also indicate how to use the Pimsner-Voiculescu six-term exact sequence 
to compute their K-theory.


\section{Notation and basic facts concerning Hecke algebras}
\label{sec:notationHA}

Suppose that $G$ is a group and
$H$ is a subgroup
of $G$. 
Then 
$(G,H)$ 
is 
a \emph{Hecke pair} 
if and only if 
$H$ 
is \emph{almost normal} 
in $G$, 
that is, 
it satisfies 
the following three 
equivalent 
conditions:
\begin{itemize}
\item 
For any 
$g\in G$, 
the index of 
$H\cap gHg^{-1}$ 
in $H$ 
and 
$gHg^{-1}$ 
is finite, 
in other words, 
$G$ equals 
the commensurator subgroup 
of $H$ 
in $G$, 
$
\{x\in G\colon x\ \text{commensurates}\ H\}$;  
\item 
Every $H$ double coset contains finitely many right $H$ cosets; 
\item 
Every $H$ double coset contains finitely many left $H$ cosets. 
\end{itemize}
That 
$(G,H)$ is 
a Hecke pair 
is sufficient 
to ensure 
that 
the product 
of any two 
$H$ double cosets 
is 
a finite union of 
$H$ double cosets. 
This condition 
also allows 
to define 
an associative algebra, 
whose 
rule of multiplication  
for 
two $H$ double cosets 
is 
a `weighted version' 
of the formula 
for  
the set-theoretic product of 
those double cosets. 

\begin{definition}[Hecke algebra over a ring]
Suppose that 
$(G,H)$ is a Hecke pair. The  \emph{Hecke algebra}  $\mathbb{Z}[G,H]$ 
of $(G,H)$ over the ring 
of integers $\mathbb{Z}$ 
is
the free $\mathbb{Z}$-module 
with basis 
the set $H\backslash G/H$ 
of double cosets 
of $H$ 
in $G$ endowed 
with the multiplication defined by
\[
HgH\cdot Hg'H=\sum_{x\in G} \alpha(g,g';x)\,HxH\,,
\]
where 
$HgH=\bigsqcup_i Hg_i$,\quad $Hg'H=\bigsqcup_j Hg'_j$ 
and 
$
\alpha(g,g';x):=\# \{(i,j)\colon g_ig'_j\in Hx\}
$. 
Also, 
for any ring $\mathcal{R}$, the
\emph{Hecke algebra of $(G,H)$ over $\mathcal{R}$}
is 
$\mathcal{R}[G,H]:=\mathbb{Z}[G,H]\otimes_\mathbb{Z}\mathcal{R}$.
See \cite[Section~4]{kri} for the details. 
When $\mathbb{Z}[G,H]$ is commutative, $(G,H)$ is called a \emph{Gelfand pair}.
\end{definition}

When 
the ring $\mathcal{R}$ 
is the field 
of complex numbers, 
Hecke algebras 
can be endowed with 
a conjugate linear involution
 *$\colon \mathbb{C}[G,H]\to  \mathbb{C}[G,H]$ 
defined on scalar multiples of 
the canonical basis 
by $(\lambda\, HxH)^*:=\overline{\lambda}\,Hx^{-1}H$ and 
turning $\mathbb{C}[G,H]$ into a unital *-algebra.  
 If the (unital, *-preserving) representations of $\mathbb{C}[G,H]$ on Hilbert space
 give rise to a C*-norm on  $\mathbb{C}[G,H]$, we define the
 universal enveloping C*-algebra $C^*(G,H)$ to be the
 completion of $\mathbb{C}[G,H]$ with respect to this norm.
 The existence of $C^*(G,H)$ is not automatic and
 we shall be interested in classifying our examples according to it.
%

%
Note that every normal subgroup of $G$ is almost normal and the Hecke algebra with respect to a normal subgroup is the group algebra of the quotient group. In this sense Hecke algebras are an alternative to group algebras of quotients and they play an analogous representation-theoretic role.

We denote 
the finite numbers 
of right 
respectively 
left cosets 
contained in 
a double coset 
by 
\[
R(HxH):=\#\,H\backslash HxH\quad 
L(HxH):=\#\,HxH/H=R(Hx^{-1}H)\,.
\]
We usually 
abuse notation 
and write 
$R(x):=R(HxH)$ 
and 
$L(x):=L(HxH)$.  
On many occasions 
we make use 
of Corollary~4.5 in~\cite{kri}, 
which states that,  
for any ring $\mathcal{R}$, 
the function $R$ 
extends to 
an $\mathcal{R}$-algebra homomorphism 
$R\colon \mathcal{R}[G,H]\to \mathcal{R}$.

\section{Notation and basic facts concerning trees and their automorphism groups}
\label{sec:notation}

Throughout the document 
$\mathcal T$ will denote 
an infinite tree 
with 
vertex set $\mathcal T^0$ 
and 
edge set $\mathcal T^1$. 
A tree 
is 
a bipartite graph, 
that is,  
we have 
a decomposition 
$\mathcal T^0= \mathcal T_{even} \cup \mathcal T_{odd}$ 
with the property that 
each vertex in one subset 
is connected only to 
vertices in the other subset. 
Vertices 
in the same subset 
are said to have 
the same \emph{type} 
and 
automorphisms of the tree 
which preserve 
this decomposition 
are said to be 
\emph{type-preserving}.

There is 
a graph-theoretic distance function $d: \mathcal T^0\times \mathcal T^0 \to \NN$ 
on the vertices 
of the tree. 
In cases 
where we use 
a distinguished vertex 
in our arguments, 
we choose 
our notation 
so that 
$\mathcal T_{even}$ is 
the set 
of vertices 
at even distance from 
our distinguished vertex 
and, 
consequently, 
$\mathcal T_{odd}$ is 
the set of vertices 
at odd distance from 
our distinguished vertex. 
Given a group 
$G$ 
of automorphisms 
of $\mathcal T$, 
the stabilizer of 
a distinguished vertex 
will be denoted 
$G_0$.

By 
a \emph{geodesic line} 
in a tree 
we will mean 
a subgraph  
isomorphic to 
the real line, 
$\mathbb{R}$, 
triangulated by 
the 
vertex set $\mathbb{Z}$ 
of integers. 
A \emph{ray} 
in a tree 
is a subgraph  
isomorphic to 
the 
half-line 
of positive reals, 
triangulated by 
the set of 
non-negative integers. 
An \emph{end} 
of a tree 
is 
an equivalence class 
of rays 
in the tree, 
where 
two rays 
are equivalent 
if their intersection 
is a ray. 

We make 
constant use 
of the classification 
of automorphisms 
of a tree 
obtained by Tits 
in Proposition~3.2 
of~\cite{tit}, 
for which 
section~6.4 in~\cite{ser-tre}
is also 
a good reference. 
An automorphism 
of a tree 
either 
(1) fixes a vertex; 
(2) inverts an edge 
(that is, 
exchanges 
its two vertices);  
or 
(3) has 
a (unique) 
invariant geodesic line, 
called its \emph{axis}, 
on which 
it induces 
a translation 
of non-zero amplitude, 
which is called 
the \emph{translation length} 
of the automorphism.  
The possibilities 
(1)--(3) above 
are 
mutually exclusive 
and 
an automorphism 
satisfying~(1) or~(2) 
is called \emph{elliptic} 
while 
an automorphim 
satisfying~(3) 
is called \emph{hyperbolic}. 
Elliptic
automorphisms 
of type~(2) 
are called 
\emph{inversions}
and the requirement that 
a group acts 
on a tree
\emph{without inversion} 
means that 
its elements act only 
as automorphisms of type~(1) or~(3). 
Of 
the two ends 
defined by 
rays 
contained in 
the axis of 
a hyperbolic automorphism, 
$h$ say,  
one 
has the property that 
every ray 
contained in 
the axis  
in that end 
is mapped 
into itself 
by $h$; 
this end
is called \emph{attracting for $h$} 
while 
the other end 
is called \emph{repelling for $h$}. 

We will be working 
with trees 
which 
turn out to be 
\emph{semi-homogenous} 
as a consequence 
of the conditions 
imposed on 
their automorphism groups. 
This means that 
vertices 
of a given type 
are incident on 
the same number 
of edges 
and 
we call this number 
the 
\emph{vertex degree} or \emph{ramification index} 
of that type 
of vertex. 
The formulae 
are more elegant 
if the ramification indices 
of 
the even 
and 
odd types 
are denoted 
$q_0+1$ and $q_1+1$ respectively. 
If $q_0=q_1$, 
the tree is \emph{homogeneous}. 
%
A group 
can only 
act transitively on 
the vertex set $\mathcal T^0$ 
if 
$\mathcal T$ is homogeneous.

We shall use 
the basic facts 
collected in 
the following Lemma 
several times 
during 
our analysis 
of Hecke algebras; 
they relate algebraic properties 
of the Hecke algebra
to geometric properties
of the group action.

\begin{lemma}\label{lem:metric-Hecke}
Suppose the group $\Gamma$ acts by isometries on the pseudo-metric space $X$ and let $o$ 
be a point of $X$ 
with stabilizer $\Gamma_0$  in $\Gamma$. Then the following statements hold for the pair $(\Gamma,\Gamma_0)$. 
\begin{enumerate}
\item\label{lem:metric-Hecke0}
The number $L(g)$ is the cardinality of the orbit of $g \udodot o$ under $\Gamma_0$. 
\item\label{lem:metric-Hecke1}
If  the sphere  with center $o$ and radius $d(o,g \udodot o)$  is finite,  then $L(g)$ is finite. 
In particular,  if  all balls in $X$ are finite, then $(\Gamma,\Gamma_0)$ is a Hecke pair. 
\item\label{lem:metric-Hecke2}
If the action of $\Gamma_0$ on the sphere with center $o$ and radius $d(o,g \udodot o)$ 
is transitive, then $\Gamma_0 g\Gamma_0$ is the set 
of elements of $\Gamma$ mapping $o$ to some element on the sphere with center $o$ 
and radius $d(o,g \udodot o)$. In particular 
$L(g)=R(g)$ is the cardinality of that sphere and 
$\Gamma_0 g\inv \Gamma_0 = \Gamma_0 g\Gamma_0$
is self adjoint in the Hecke algebra 
$\C[\Gamma,\Gamma_0]$.
\item\label{lem:metric-Hecke3}
If all balls in $X$ are finite and the action of $\Gamma_0$ on all spheres with center $o$ 
is transitive, then $(\Gamma,\Gamma_0)$ is a Gelfand pair, i.e. $\C[\Gamma,\Gamma_0]$
is commutative. 
\end{enumerate}
\end{lemma}
\begin{proof}
First note that for all $g_1,g_2\in\Gamma$ 
\[
g_1\udodot o =g_2\udodot o \iff g_2^{-1}g_1\udodot o = o 
\iff  g_2^{-1}g_1\in\Gamma_0 \iff  g_1\in g_2\Gamma_0.
\]
Thus each left $\Gamma_0$-coset $g\Gamma_0$ corresponds to an element $g\udodot o$  in the $\Gamma$-orbit of $o$. Now~(\ref{lem:metric-Hecke0}) follows from the observation that $\Gamma_0g\Gamma_0=\Gamma_0(g\Gamma_0)$.

Since the action is by isometries on $X$, the  orbit 
of $g \udodot o$ under $\Gamma_0$ is contained in  
the sphere with center $o$  and radius $d(o,g \udodot o)$, 
hence (\ref{lem:metric-Hecke1}) follows from  (\ref{lem:metric-Hecke0}).
We now turn to~(\ref{lem:metric-Hecke2}).  
Since the action of $\Gamma_0$ on the
 sphere with center $o$ and radius $d(o,g \udodot o)$ is transitive, $\Gamma_0 g\Gamma_0$
is the set of elements of $\Gamma$ mapping $o$ to some element 
on this sphere. 
In particular $d(g\inv \udodot o, o) = d(o,g \udodot o)$, so $g\inv $ is such an element and
 $\Gamma_0 g\inv\Gamma_0 =\Gamma_0 g \Gamma_0$.
This also implies that  $L(g)$ and $ R(g)$ are both equal to the cardinality of the sphere. 
Statement ~(\ref{lem:metric-Hecke3}) now follows because the generators of $\mathbb{C}[\Gamma,\Gamma_0]$
are self-adjoint by (\ref{lem:metric-Hecke2}), and thus satisfy $a b = a^* b^* = (ba)^* = ba$.
\end{proof}

\begin{remark}\label{rem:metric-Hecke}
In Lemma~\ref{lem:metric-Hecke} we may assume that 
 that $\Gamma$ acts transitively
 by taking 
 $X=\Gamma \udodot o$.  
\end{remark} 

\begin{remark}\label{rem:Hecke_subgroup=Stab(point)}
Even in cases when the action of $\Gamma _0$ 
on spheres centered at $o$ is not transitive, the triangle inequality imposes
a restriction on the scattering resulting from the multiplication of double cosets. Indeed, 
suppose that $\Gamma$ acts by isometries on a metric space $X$ in which spheres  are finite, and 
let $A_n$ be the $\C$-linear span of all the double $\Gamma_0$-cosets  of
the elements $\gamma \in \Gamma$ 
such that $d(\gamma \udodot o,o)= n$ (with the convention that $A_n=\{0\}$ if no such elements exists). 
Then  $A_n A_m\subseteq \bigoplus_{i=|m-n|}^{m+n}A_{i}$  in $\mathbb{C}[\Gamma,\Gamma_0]$. 
\end{remark}

\section{Edge-transitive automorphism groups 
relative to a vertex stabilizer}
\label{sec:Hecke_relStab(vertex)}
In this section we will assume that the locally finite tree  $\mathcal{T}$
has at least three edges at each vertex, and we select a distinguished vertex $o$ of $\mathcal T$.
We will also suppose that $\Gamma$ is a group that acts on $\mathcal T$ in such a way that the
stabilizer $\Gamma_0$ of $o$ is a proper subgroup of $\Gamma$ that acts transitively on 
each sphere about $o$. 
We will show that under these conditions 
the underlying tree $\mathcal{T}$ is necessarily 
semi-homogenous and that such a pair $(\Gamma, \Gamma_0)$ is a Gelfand pair
for which the Hecke algebra is simply the algebra of polynomials in a self adjoint 
variable.
The situation under consideration includes as special cases those
in which $\Gamma$ is  the full group of automorphisms
of a homogeneous tree,  the subgroup of type-preserving 
automorphisms of a semi-homogenous tree, and 
many of their proper subgroups, some of which cases
have been studied before independently e.g.~in \cite{TAct.harm.c}.
More examples of the present situation arise
naturally from  semisimple matrix groups  
of rank~$1$
over local fields acting on their Bruhat-Tits trees, 
the simplest example 
being $\mathrm{SL}_2(\mathbb{Q}_p)$, see Remark~ \ref{rm:endofsection3}. 

\begin{proposition}\label{pro:3edges}
Let $\mathcal T$ be a tree 
with at least $3$ edges 
at each vertex 
and let
$o$ be 
a distinguished vertex
of $\mathcal T$.
Suppose that  
$\Gamma$ is a group 
acting by automorphisms 
of $\mathcal T$  
and assume that 
the stabilizer $\Gamma_0$ 
of $o$ is a proper subgroup of $\Gamma$ that
acts transitively 
on each sphere 
around $o$.
Then 
the tree $\mathcal T$ 
is semi-homogeneous  
and
the action of $\Gamma$ on 
$\mathcal T^0$ 
is either transitive or else 
has two orbits, namely,  the sets 
$\mathcal{T}_{even}=\{v\in \mathcal{T}^0\colon d(v,o)\ \text{is even}\}$ 
and 
$\mathcal{T}_{odd}=\{v\in \mathcal{T}^0\colon d(v,o)\ \text{is odd}\}$. 
\end{proposition}

\begin{proof}  
By
part~(4)  of Lemma~\ref{lem:metric-Hecke}  
applied to 
the action 
of $\Gamma$  
on the subset $\mathcal{T}^0$,  
$(\Gamma,\Gamma_0)$ 
is a Gelfand pair. 
In order 
to verify 
that 
the tree $\mathcal{T}$ 
is semi-homogenous 
and that 
the action 
of $\Gamma$ 
on the set $\mathcal{T}^0$ 
of vertices 
of $\mathcal{T}$ 
has 
at most~$2$ orbits, 
it is 
enough 
to show 
the following lemma. This will finish the proof of the proposition.
\end{proof}
\begin{lemma}\label{lem:Gamma transports by even distances}
For each even integer $D$ 
  there exists  an elliptic element $\delta\in \Gamma$ 
with $d(\delta \udodot o,o)=D$. 
\end{lemma}

\begin{proof}[Proof of Lemma~\ref{lem:Gamma transports by even distances}]
By  \cite[Proposition~3.4]{tit} 
$\Gamma$ contains a hyperbolic element $\eta$. 
In particular, $d(\eta^n \udodot o,o)$ goes to infinity 
as $n$ goes to infinity, and one can choose
 an element $\gamma$ in $\Gamma$ such that 
$d(\gamma \udodot o,o)\ge D$. 
Let $v$ be the vertex 
at distance $D/2$ from~$\gamma \udodot o$ 
on the path 
joining~$o$ to~$\gamma \udodot o$.
Let  $w$ be a vertex  
at distance $D/2$ from $v$ 
on the sphere with radius $d(\gamma \udodot o,o)$ 
such that $[\gamma \udodot o,o]\cap [o,w]\cap [w,\gamma \udodot o]=\{v\}$. 
Such a vertex~$w$ exists 
because there are 
at least $3$ edges at each vertex 
by assumption.
The situation is illustrated in the following diagram; note that 
because of our choice of $w$ we have $d(w,\gamma \udodot o) = D$. 

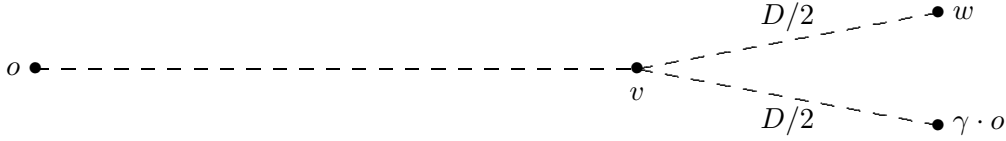
\centerline{
\beginpicture
\setcoordinatesystem units <1cm,0.5cm>  
\setplotarea x from -9 to 5, y from -2.5 to 2.5            
\put{$\bullet$}  at  0 0
\put{$\bullet$}  at -8 0
\put{$\bullet$}  at  4 1.5
\put{$\bullet$}  at  4 -1.5
\put{$o$}  [r]   at -8.2 0
\put{$\gamma\udodot o$}  [l]   at 4.2 -1.5
\put{$w$}  [l]   at 4.2 1.5
\put{$v$}  [t]   at 0 -0.5
\put{$D/2$}  [t]   at 2 -1
\put{$D/2$}  [b]   at 2 1
\setdashes
\putrule from -8 0 to 0 0
\plot 4 1.5   0 0  /
\plot 4 -1.5   0 0  /
\endpicture
}


Choose an element $\pi$  of $\Gamma_0$ mapping $\gamma \udodot o$ to $w$; 
such an element exists because  $\Gamma_0$ acts transitively on the sphere 
with center $o$ and radius $d(\gamma \udodot o,o)$. 
The element $\delta:=\gamma^{-1}\pi\gamma$ 
is elliptic,  because $\pi$ is, and since 
\[
d(\delta\udodot o, o) =  d( \gamma^{-1}\pi\gamma\udodot o, o) = 
d(\gamma^{-1}  \udodot w, o) =  d(w, \gamma \udodot o) =D,
\]
this finishes the proof of the lemma and of \proref{pro:3edges}
\end{proof} 

\begin{theorem}\label{thm:rel-vertexStab}
Under the same assumptions as in Proposition~\ref{pro:3edges},
let $q_0+1$ (respectively $q_1+1$) denote 
the degrees 
of the even (respectively odd) vertices 
and 
put $q:=q_0$  if $\mathcal{T}$ is homogenous.  
Define $\Gamma_n := \{ \gamma \in \Gamma: d(\gamma \udodot o , o) = n\}$.  
Then
\begin{enumerate} 
\item in the case when $\Gamma$ acts transitively  on vertices,
 the canonical basis 
of $\mathbb{C}[\Gamma,\Gamma_0]$ 
is the set $\{\Gamma_n\colon n\in\mathbb{N}\}$,
the map $T \mapsto \Gamma_1$ 
 extends to an isomorphism of the algebra $\mathbb{C} [T]$  of complex 
 polynomials in the variable $T$ onto 
 $\mathbb{C}[\Gamma,\Gamma_0]$, and 
 the elements 
of the canonical basis 
of $\mathbb{C}[\Gamma,\Gamma_0]$ 
satisfy the recursion relations
\[
\Gamma_1\Gamma_n = 
(q+\delta_{n,1})\,\Gamma_{n-1} + \Gamma_{n+1} \qquad \text{ for } n\ge 1;
\] 

\item in the case when $\Gamma$ has two orbits of vertices,
 the canonical basis 
is the set $\{\Gamma_{2n}\colon n\in\mathbb{N}\}$, the isomorphism of 
$\mathbb{C} [T]$ to   $\mathbb{C}[\Gamma,\Gamma_0]$
is determined by $T \mapsto \Gamma_2$ and the recursion
relations are   
\[
\Gamma_2\Gamma_{2n}=
(q_0+\delta_{n,1})q_1\,\Gamma_{2n-2}+(q_1-1)\,\Gamma_{2n}+\Gamma_{2n+2} \qquad \text{ for }
 n\ge 1
.\]
\end{enumerate}
\end{theorem}
\begin{proof}
That the given sets are canonical bases is immediate from 
parts~(3) and (4) of Lemma~\ref{lem:metric-Hecke} 
applied to the action 
of $\Gamma$  
on the set $\mathcal{T}^0$ of vertices.
To show that $T\mapsto \Gamma_1$ and $T\mapsto \Gamma_2$ 
give isomorphisms we need the following lemma.

\begin{lemma}\label{lem:reach-index_sum}
If both $\Gamma_n$ and $\Gamma_m$ 
are non-empty, 
then 
there exist $\gamma_n\in\Gamma_n$ and $\gamma_m\in\Gamma_m$ 
such that $\gamma_n\gamma_m\in\Gamma_{n+m}$.
\end{lemma}
\begin{proof}[Proof of Lemma~\ref{lem:reach-index_sum}]
%
Considering the inverse of the product $\gamma_n\gamma_m$ if necessary, 
we may suppose that $n\leq m$. 
Choose $\gamma_m$ in $\Gamma_m$ arbitrarily. Suppose first that $\Gamma_n$
has an elliptic element, say $\delta_n$.
 Then by transitivity of the action of $\Gamma_0$ on the sphere 
of radius $n$ around $o$, a suitable conjugate $\gamma_n$ of $\delta_n$ 
by an element of $\Gamma_0$ satisfies $\gamma_n \udodot o\notin [o,\gamma_m]$. 
With this choice of $\gamma_n$ we have $d(\gamma_n\gamma_m \udodot o,o)=n+m$. 
Alternatively, suppose $\Gamma_n$  has a hyperbolic element, say $\delta_n$.
 Then, again by transitivity of the action of $\Gamma_0$ on the 
sphere of radius $n$ around $o$, there is a $\Gamma_0$-conjugate, 
$\gamma_n$ say, of $\delta_n$ such that the closest point on its axis to $o$ 
lies on $[o,\gamma_m \udodot o]$ 
and 
$o$ does not lie 
between $\gamma_n.o$ and $\gamma_m.o$.
In this case, this choice of $\gamma_n$ guarantees 
$d(\gamma_n\gamma_m \udodot o,o)=n+m$, finishing the  proof of \lemref{lem:reach-index_sum}.
\end{proof}

We continue 
the proof 
of Theorem~\ref{thm:rel-vertexStab}.
%
Let $n$ and $m$ be integers with $n\leq m$. By Lemma~\ref{lem:reach-index_sum} 
and Remark~\ref{rem:Hecke_subgroup=Stab(point)} we have 
\begin{equation}\label{eq:justderived}
\Gamma_n\Gamma_m=\sum_{i=m-n}^{m+n} \lambda_i\,\Gamma_i;\quad \lambda_{m+n}\neq 0\,.
\end{equation}
If $\Gamma$ is transitive on the set of vertices of $\mathcal{T}$, we use this equation with $n=1$ 
and arbitrary $m\ge 1$. If $\Gamma$ has two orbits of vertices, then $\Gamma_m$ is empty for $m$ odd  and we use this equation with $n=2$ and $m$ even and at least $2$. 

Since 
$\mathbb{C}[\Gamma,\Gamma_0]$ is abelian,
\eqref{eq:justderived} implies that
each element 
of $\mathbb{C}[\Gamma,\Gamma_0]$ 
is a polynomial in $\Gamma_i$, with $i=0$ 
if $\Gamma$ acts transitively 
on the vertex set $\mathcal T^0$ 
and  $i=2$
if the action of $\Gamma$ 
on $\mathcal T^0$ 
has two orbits. This  implies that the map 
$T \mapsto \Gamma_i$ extends to an algebra homomorphism
of the algebra $\mathbb{C}[T]$ of polynomials on $T$ onto $\mathbb{C}[\Gamma,\Gamma_0]$.
Since every proper quotient of $\mathbb{C}[T]$
is a finite-dimensional $\mathbb{C}$-vector space, 
and both canonical bases are infinite
we have an isomorphism in both cases.

Next we verify
the recursion relations. 
If 
$\Gamma$ acts 
transitively on $\mathcal{T}^0$, then  
the set-theoretic product 
of the double cosets 
$\Gamma_1$ and $\Gamma_n$  
is 
 $\Gamma_{n-1}\cup\Gamma_{n+1}$.
Then 
\[
\Gamma_1\Gamma_n= a\Gamma_{n-1}+b\Gamma_{n+1}
\] 
with both  $a$ and~$b$ positive integers. 
On applying the algebra homomorphism $R$ 
to this equation, noting that $R(\Gamma_0)=1$ 
and $R(\Gamma_n)=(q+1)q^{n-1}$ 
for $n\ge 1$ by  part~(3) of Lemma~\ref{lem:metric-Hecke}, 
we obtain the following equations 
for the cases $n=1$ and~$n>1$: 
\[
(q+1)^2q^{n-1}=
\begin{cases}
a+b(q+1)q& \text{if}\ n=1\,;\\
a(q+1)q^{n-2}+b(q+1)q^n& \text{if}\ n>1\,.
\end{cases}
\]
In these equations 
we cannot have $b\ge 2$ 
because then 
$b(q+1)q^n\ge 2(q+1)q^n$, 
which is 
strictly larger than the left hand side $(q+1)^2q^{n-1}$, 
whatever the value of $n$; thus we conclude that $b=1$. 
This enables us 
to solve 
the equations 
for the parameter~$a$
in both cases, from which we obtain
$a = q+1 $ when $n = 1$ and $a = q$ when $n >1$. 

If $\Gamma$ acts with two orbits on $\mathcal{T}^0$, then
the set-theoretic product of the double cosets 
$\Gamma_2$ and $\Gamma_{2n}$ 
is the set $\Gamma_{2n-2}\cup\Gamma_{2n}\cup\Gamma_{2n+2}$. 
Then 
\[
\Gamma_2\Gamma_{2n}= a\Gamma_{2(n-1)}+b\Gamma_{2n}+c\Gamma_{2(n+1)}
\]
with $a$, $b$, and $c$ positive integers.
Applying again the homomorphism $R$ using 
 part~(3)  of Lemma~\ref{lem:metric-Hecke} 
 and the cardinality of the spheres of a semi-homogeneous tree
 in this case, we see that 
\[
(q_0+1)^2q_1^2 
= a  + b (q_0+1) q_1 + c (q_0+1)q_0 q_1^2
\qquad \text{ for \ } n = 1
\]
while
\[
(q_0+1)q_1 (q_0+1) \prod_{j = 1}^{2n-1} q_j
= a (q_0+1) \prod_{j = 1}^{2n-3} q_j + b (q_0+1) \prod_{j = 1}^{2n-1} q_j + c (q_0+1) \prod_{j = 1}^{2n+1} q_j
\]
for $n>1$.
Arguing as above, 
we conclude  that $c=1$ for every $n$. 

In order to determine~$b$
we argue as follows.
The product 
of two representatives 
$\gamma_2$ in $\Gamma_2$ 
and~$\gamma_{2n}$ in $\Gamma_{2n}$ 
lies in $\Gamma_{2n}$, 
if and only if 
$\gamma_{2n}$ 
moves 
the vertex~$o$ 
by $2n$ units,  
$\gamma_{2}$ 
moves~$o$   
by $2$ units, 
and 
their product 
moves~$o$ 
by $2n$ units. 
This happens 
if and only if 
the segments $[o,\gamma_{2n}\udodot o] $ and~$[o,\gamma_2\inv \udodot o]$ 
intersect in precisely one edge, which gives us 
$b=(q_1-1)$ possibilities for the representative $\gamma_2$, 
as claimed in statement~(3)(b) of our theorem. 
Finally $a = (q_0 +\delta_{n,1})q_1$ follows from the above equation.
\end{proof}

The multiplication table 
for the canonical basis 
of $\mathbb{C}[\Gamma,\Gamma_0]$ 
may be derived 
from the recursion relations 
by induction 
on the parameter~$n$. 
We leave the details 
to the reader. 
\begin{corollary}
The complete multiplication table 
for the elements 
of the canonical basis 
of $\mathbb{C}[\Gamma,\Gamma_0]$ 
is given by: 
\begin{enumerate}
\setcounter{enumii}{2}
\item if 
$\Gamma$ acts transitively 
on $\mathcal{T}^0$, then
\[
\Gamma_n\Gamma_m = \Gamma_m\Gamma_n=
\Gamma_{m+n}+q^{n-1}(q+\delta_{mn})\Gamma_{m-n}+(q-1)\sum_{l=1}^{n-1} q^{l-1} \Gamma_{m+n-2l}
\]
for $m\ge n>0$; 
\item if 
$\Gamma$ acts 
with two orbits 
on $\mathcal{T}^0$,
letting
$q_l = q_0$ 
for $l$ even 
and 
$q_l = q_1$ 
for $l$ odd, 
then 
\[
\Gamma_{2n}\Gamma_{2m} = \Gamma_{2m}\Gamma_{2n}=
\Gamma_{2(m+n)}+
q_1^nq_0^{n-1}(q_0+\delta_{mn})\Gamma_{2(m-n)}
+\sum_{l=1}^{2n-1} (q_{l}-1)\prod_{i=1}^{l-1} q_{i} \Gamma_{2(m+n-l)}
\]
for $m\ge n >0$.
\end{enumerate}
\end{corollary}

We close 
this section 
with several remarks. 

\begin{remark}\label{rm:endofsection3}
~\\[-3ex]
\begin{enumerate}
\item 
It follows from 
part~(1) 
of Theorem~\ref{thm:rel-vertexStab}  
that $\mathbb{C}[\Gamma,\Gamma_0]$ 
does not have 
an enveloping $C$*-algebra 
because 
evaluation 
of polynomials 
at real numbers 
gives a set 
of *-representations 
of $\mathbb{C}[\Gamma,\Gamma_0]\cong\mathbb{C}[T]$, 
whose norms 
are not uniformly bounded. 
\item 
It is often convenient 
to work with 
`normalized double cosets', 
which are obtained by 
scaling 
a double coset 
with the inverse 
of the value   
that the $R$-function 
attains there, 
thus putting 
$\tilde \Gamma_k := R(\Gamma_k)\inv \Gamma_k$. 
The recurrence relations 
for these normalized 
double cosets 
become 
\begin{align*}\label{eqn:fullhomnormalizedreln}
\tilde\Gamma_1 \tilde\Gamma_n &= 
\frac{1}{q+1}\biggl(\tilde\Gamma_{n-1} + q \tilde\Gamma_{n+1}\biggr), \quad n\geq 1\,;\\
\tilde\Gamma_2\tilde\Gamma_{2n} &=
\frac{1}{(q_0+1)q_1}
\biggl(
\tilde\Gamma_{2n-2}+(q_1-1)\tilde\Gamma_{2n}+ q_0q_1\tilde\Gamma_{2n+2}
\biggr)
; \quad n\ge 1\,.
\end{align*}
The first set 
of recurrence relations 
are those 
obtained for 
the spherical functions 
on the homogeneous tree $\mathcal T$ 
by Fig{\`a}-Talamanca and Nebbia 
in~\cite{harmAna+repThry(Gs-act>homogTrees)}. 
\item 
Theorem~\ref{thm:rel-vertexStab} 
applies to 
the action 
of $\mathrm{SL}_2(\mathbb{Q}_p)$ 
on its Bruhat-Tits tree. 
For this action 
there is a vertex 
whose stabilizer 
is the subgroup  $\mathrm{SL}_2(\mathbb{Z}_p)$. 
That 
the Hecke algebra 
of the pair 
$(\mathrm{SL}_2(\mathbb{Q}_p),\mathrm{SL}_2(\mathbb{Z}_p))$ 
does not have a universal $C$*-algebra
was observed by Hall 
in Proposition~2.21 
of her PhD thesis~\cite{hall}.
\end{enumerate}
\end{remark}

\section{Strongly transitive automorphism groups
relative to an edge fixator.}
\label{sec:Hecke_relFix(edge)}
In this section we assume that  $\mathcal{T}$ is a locally 
finite tree with a distinguished edge $e$.
The theme of the section is the structure 
of the Hecke algebra of a given group of automorphisms of  $\mathcal{T}$ 
with respect to the subgroup that fixes $e$.  Motivated by the case of
algebraic groups, shall impose an extra transitivity assumption on the  
type-preserving automorphisms of the given group.

\begin{definition}
A type-preserving group of automorphisms of a tree is \emph{strongly transitive\/} 
if it acts transitively on the set of 
doubly-infinite geodesics in the tree 
and if the stabilizer of some doubly-infinite geodesic 
$\mathcal{A}$ acts transitively
on the set of edges in $\mathcal{A}$. 
\end{definition}

Although the current context includes groups that are not type-preserving, we have chosen to restrict the definition of strong transitivity to actions of type-preserving automorphisms,
for consistency with the established literature.
This restriction is natural in the context of algebraic groups because of the link with $BN$-pairs. In that context, releasing the type-preserving assumption would mean that the Weyl group would no longer be a Coxeter group in general with significant consequences for the last $BN$-pair condition. 
It will be convenient to list here some
 reformulations of the property of strong transitivity. The proof is left as an exercise for the reader.

\begin{lemma}\label{lem:strongly-transitive}
Suppose that 
$G$ is 
a group of 
type preserving automorphisms of 
a locally finite tree $\mathcal T$; 
then 
the following conditions 
are equivalent.
\begin{enumerate}
\item [(1)] 
$G$ is strongly transitive.
\item [(2)] $G$ acts transitively on the edge set $\mathcal T^1$ 
and there exists an edge $f\in \mathcal T^1$ whose stabilizer 
acts transitively on the set of all doubly-infinite geodesics through~$f$. 
\item [(3)] 
$G$ acts transitively on the edge set $\mathcal T^1$ 
and  the stabilizer of each edge acts transitively on the set of all 
doubly-infinite geodesics through that edge. 
\item [(4)] 
For any two doubly-infinite geodesics $\mathcal A_1$ and $\mathcal A_2$ 
and any two edges $f_1\in\mathcal A_1$ and $f_2\in\mathcal A_2$, 
there exists an element $g\in G$ such that $g \cdot \mathcal A_1 = A_2$
and $g \cdot f_1=f_2$. 
\end{enumerate}
\end{lemma}

 It follows from \lemref{lem:strongly-transitive} 
 that the stabilizer of an arbitrary doubly-infinite geodesic is transitive on its set of edges.

Next we consider group automorphisms of a locally  finite tree such that 
the subgroup of type-preserving automorphisms is strongly transitive.
We will show that  the
stabilizer subgroup of the distinguished edge $e$ is a Hecke subgroup
and we will characterize the corresponding Hecke algebra.
This  context includes,  for instance,
the group of all automorphisms (respectively all type-preserving automorphisms) 
of a (semi-)homogenous tree and also natural examples 
arising from semisimple matrix groups of rank~$1$ over local fields, 
for which the underlying trees are the associated Bruhat-Tits trees.
Specifically, working along the lines of  \cite{IHES_25_5_0},
we will show that 
when the subgroup 
$G^+$ of  type preserving automorphisms is strongly transitive,
$(G, \mathcal{B})$ and $(G^+, \mathcal{B})$ are Hecke pairs and 
we will derive the structure of  their Hecke algebras.  

\begin{proposition}\label{pro:rel-edgeStab}
Let $G$ be a group acting by 
 automorphisms of a tree $\mathcal{T}$ 
 in such a way that the subgroup $G^+$ 
of type-preserving automorphisms is 
strongly transitive.
For a distinguished edge $e$ of $\mathcal{T}$, 
let $\mathcal{B}$ be the subgroup of $G$ 
fixing $e$ pointwise and let $\mathcal A$ be a doubly-infinite geodesic
containing  $e$. 

Then $\mathcal{T}$ is semi-homogenous and the action of $G$ 
on $\mathcal{T}^0$ is either transitive (iff $G\neq G^+$) or else 
has two orbits (iff $G= G^+$).
If $G\neq G^+$ there is in $G$ an inversion $i$
of the edge $e$ that stabilizes $\mathcal{A}$, yielding  a short exact sequence
$
1\longrightarrow G^+\longrightarrow G\longrightarrow \mathbb{Z}/2\mathbb{Z}\longrightarrow 1\,.
$
This short exact sequence splits if and only if 
$G$ contains an inversion that is 
an involution. 
Furthermore, both $(G,\mathcal B)$ and $(G^+,\mathcal B)$ are Hecke pairs.
\end{proposition}

\begin{proof}
If $G=G^+$ everything except the statement about Hecke pairs follows easily from 
\lemref{lem:strongly-transitive}. 
Suppose now that $G\neq G^+$. Then there exists $g\in G$ which is not type-preserving.  
Using \lemref{lem:strongly-transitive} (4) we can then construct an element $i\in G$ which inverts the edge $e\in\mathcal T^1$
and  which stabilizes the doubly-infinite geodesic $\mathcal{A}$ containing $e$. %
Then 
\begin{equation}\label{rem:G=G+N(B)}
G= G^+\cup G^+i
\end{equation}
from which it follows that there exists a short exact sequence
\[
1\longrightarrow G^+\longrightarrow G\longrightarrow \mathbb{Z}/2\mathbb{Z}\longrightarrow 1\,
\]
that splits if and only if $G$ contains an inversion that is an involution.

Next we verify that the pairs of groups 
$(G^+,\mathcal{B})$ and $(G,\mathcal{B})$ are both Hecke pairs.
Let $X$ be the metric space whose set of points 
consists of  the midpoints of all edges 
of $\mathcal{T}$ (defined combinatorially as 
the set of unordered pairs $\{f,\overline{f}\}$, 
where $f$ runs through all edges of $\mathcal{T}$)
and whose distance function 
assigns the number $n$ to the pair 
consisting of the midpoint of $f_1$ and the midpoint of $f_2$ 
whenever the longest geodesic path 
joining some vertex of $f_1$  to 
some vertex of $f_2$ has length $n+1$. 

The group $G$, and hence also 
its subgroup $G^+$,  act by isometries on $X$. 
Denote the midpoint of the edge $e$ by $o$. 
The stabilizer of $o$  in $G^+$ (respectively in $G$  )
is 
$\mathcal{B}$ 
(respectively 
the normalizer $N_G(\mathcal{B})$ of $\mathcal B$ in $G$).
Since $\mathcal T$ is locally finite, 
Part~(\ref{lem:metric-Hecke1})  of Lemma~\ref{lem:metric-Hecke} implies that 
$(G^+,\mathcal{B})$ and $(G,N_G(\mathcal{B}))$ are Hecke pairs. 
Since the index of $\mathcal{B}$ in $N_G(\mathcal{B})$ 
is at most $2$, it follows that 
$(G,\mathcal{B})$ is a Hecke pair also in the case $G \neq G^+$.
\end{proof}

\begin{theorem}\label{thm:rel-edgeStab}
Under the same conditions of hypothesis of Proposition~\ref{pro:rel-edgeStab}
we have the following:
\begin{enumerate}
\item 
The algebras $\mathbb{C}[G,\mathcal{B}]$ and  $\mathbb{C}[G^+,\mathcal{B}]$ 
are isomorphic to 
the group algebra 
of the infinite dihedral group, $D_\infty$.
\item 
The canonical basis 
of the Hecke algebras 
$\mathbb{C}[G^+,\mathcal{B}]$ 
(respectively 
$\mathbb{C}[G,\mathcal{B}]$ if $G\neq G^+$) 
consisting of 
the double cosets 
with respect to $\mathcal{B}$,  
are 
parametrized by 
the automorphism group 
of the bipartite graph (
respectively 
the abstract graph)
underlying
$\mathcal{A}$. 
More precisely: 
\begin{enumerate} 
\item 
Denote 
the stabilizer 
of the geodesic $\mathcal{A}$ 
in $G^+$ 
by $N^+$ 
and 
the fixator of $\mathcal{A}$ 
in $G^+$ 
by $T$. 
Then 
\begin{enumerate}
\item 
every automorphism 
of  $\mathcal{A}$ 
as a bipartite graph 
is realized by 
some element of $N^+$ 
and hence 
$N^+/T$ is isomorphic to 
the group of automorphisms 
of the bipartite graph $\mathcal{A}$, 
the infinite dihedral group $D_\infty$; 
\item 
we have $T=N^+\cap\mathcal{B}$, 
and therefore, 
given an element $n$ of $N^+$ 
with image $w$ 
in the quotient group $W:=N^+/T$,  
the left, right and double cosets
of all elements in $nT$ 
with respect to $\mathcal{B}$ 
are equal 
and 
we may consequently write them 
as $w\mathcal{B}$, $\mathcal{B}w$ and $\mathcal{B}w\mathcal{B}$ 
respectively; 
\item 
with the conventions 
justified in~(\romannumeral2) above, 
the group $G^+$ 
is the disjoint union $\bigsqcup_{w\in W} \mathcal{B}w\mathcal{B}$. 
\end{enumerate}
\item 
If $G\neq G^+$, 
denote 
the stabilizer 
of the geodesic $\mathcal{A}$ 
in $G$ 
by $N$ 
and 
the fixator 
of $\mathcal{A}$ in $G$ 
by $T$. 
Then 
\begin{enumerate}
\item 
every graph automorphism 
of  $\mathcal{A}$ 
is realized by 
some element of $N$ 
and 
$\widetilde{W}:=N/T$ is isomorphic to 
the group of automorphisms 
of the graph $\mathcal{A}$, 
which is 
the split extension 
of $N^+/T$ 
by the group 
$N/N^+\cong \mathbb{Z}/2\mathbb{Z}$, 
and is 
abstractly isomorphic to 
the infinite dihedral group $D_\infty$; 
\item 
since $T$ is contained in $\mathcal{B}$, 
given an element $\widetilde{n}$ of $N$ 
with image $\widetilde{w}$ 
in the quotient group $\widetilde{W}:=N/T$,  
the left, right and double cosets
of all elements in $\widetilde{n}T$ 
with respect to $\mathcal{B}$ 
are equal 
and 
we may consequently write them 
as $\widetilde{w}\mathcal{B}$, $\mathcal{B}\widetilde{w}$ and $\mathcal{B}\widetilde{w}\mathcal{B}$ 
respectively; 
\item 
with the conventions 
justified in~(\romannumeral2) above, 
the group $G$ 
is the disjoint union $\bigsqcup_{\widetilde{w}\in \widetilde{W}} \mathcal{B}\widetilde{w}\mathcal{B}$. 
\end{enumerate}
Denote 
by~$s$ and~$t$ 
the generators 
of $D_\infty$ 
that map 
the edge $e$ 
to its two neighbors 
in~$\mathcal{A}$, 
the elements of 
$\mathbb{C}[G^+,\mathcal{B}]$ 
and 
$\mathbb{C}[G,\mathcal{B}]$ 
defined by $\mathcal{B}w\mathcal{B}$ 
for $w$ in $W$ and $\widetilde{W}$ 
by $\Delta_w$ 
and, 
if $G\neq G^+$  
the element 
of $\mathbb{C}[G,\mathcal{B}]$ 
with underlying set 
$\mathcal{B}i\mathcal{B}=\mathcal{B}i=i\mathcal{B}$ 
by $\Delta_i$. 
The image of 
the element $i$ 
in the group $\widetilde{W}$ 
will also 
be denoted by~$i$; 
it 
interchanges 
the odd and even vertex 
of the edge $e$ 
and 
conjugates $s$ to $t$. 

The algebra $\mathbb{C}[G^+,\mathcal{B}]$ 
is generated 
by the identity, 
$\mathcal{B}$, $\Delta_s$ and~$\Delta_t$, 
and 
in case $G\neq G^+$,  
the algebra $\mathbb{C}[G,\mathcal{B}]$ 
is generated 
by the identity, 
$\mathcal{B}$, $\Delta_i$ and~$\Delta_s$. 
More precisely: 
Suppose that 
$w$ is an element of $W$. 
If $s_1\cdots s_n$ is 
a reduced decomposition 
of $w\in W$ 
as a word in $\{s,t\}$  
then 
$\Delta_w=\Delta_{s_1}\cdots \Delta_{s_n}$ 
and 
$R(\Delta_w)=\prod_{i=1}^n R(\Delta_{s_i})$. 
Furthermore, 
for every $\widetilde{w}$ in $\widetilde{W}$, 
we have 
$\Delta_i\Delta_{\widetilde{w}}=\Delta_{i\widetilde{w}}$.  
\end{enumerate}
\item 
\begin{enumerate}
\item 
Left and right multiplication 
with the generators $\Delta_r$; $r\in \{s,t\}$ 
on the basis 
of the $\mathbb{Z}$-module underlying $\mathbb{Z}[G^+,\mathcal{B}]$ 
that consists of 
the elements $\Delta_w$; $w\in W$ 
is given by 
\begin{align*}
\rule{0pt}{0pt}\qquad
\Delta_r\Delta_w &= q_r\Delta_{rw}+(q_r-1)\Delta_w &
\text{if 
$w$ begins with r}\,,\\
\rule{0pt}{0pt}\qquad
\Delta_w\Delta_r &= q_r\Delta_{wr}+(q_r-1)\Delta_w &
\text{if 
$w$ ends with r}\,,\\
\rule{0pt}{0pt}\qquad
\Delta_r\Delta_w &=\Delta_{rw} &
\text{if 
$w$ doesn't begin with r}\,,\\
\rule{0pt}{0pt}\qquad
\Delta_w\Delta_r &=\Delta_{wr} &
\text{if 
$w$ doesn't end with r}\,.
\end{align*}
Furthermore, 
setting $q_r:=R(\Delta_r)$, 
the set of relations 
\begin{align*}
&\Delta_r^2-\bigl(q_r+(q_r-1)\Delta_r\bigr) &(r\in \{s,t\})
\end{align*}
define a presentation of $\mathbb{Z}[G^+,\mathcal{B}]$ 
as a ring with unit. 
\item 
Suppose that $G\neq G^+$. 
For $w$ in $W$ 
write $\overline{w}$ 
for the image of $w$ 
under the automorphism 
of $W$ 
which exchanges 
the generators $s$ and $t$; 
this automorphism 
is implemented by 
conjugation with 
the element $iT$ 
in the group $\widetilde{W}$. 
The group $\mathbb{Z}/2\mathbb{Z}$ 
acts on $\mathbb{Z}[G^+,\mathcal{B}]$ 
with the non-trivial element 
therein 
sending $\Delta_w$ to $\Delta_{\overline{w}}$ 
for $w$ in $W$. 
The Hecke algebra $\mathbb{Z}[G,\mathcal{B}]$ 
is isomorphic to 
the twisted tensor product 
of this action 
of $\mathbb{Z}/2\mathbb{Z}$ 
on $\mathbb{Z}[G^+,\mathcal{B}]$. 
Furthermore, 
the set of relations 
\begin{align*}
&\Delta_r^2-\bigl(q_r+(q_r-1)\Delta_r\bigr) &(r\in \{s,t\})\,,\\
&\Delta_i^2-1\,\\
 &\Delta_i\Delta_s\Delta_i^{-1}-\Delta_{t}
\end{align*}
define a presentation of $\mathbb{Z}[G,\mathcal{B}]$ 
as a ring with unit. 
\end{enumerate}
\end{enumerate}
\end{theorem}

\begin{proof}
We begin by verifying the statements 
about  the canonical bases 
of the Hecke algebras $\mathbb{C}[G^+,\mathcal{B}]$ and  $\mathbb{C}[G,\mathcal{B}]$ 
made in claim~(2), 
starting with part~(a) 
of that claim. 

Since  
the group $G^+$ 
acts strongly transitively on $\mathcal{T}$,  
the group $N^+$ 
is transitive 
on the set of edges 
of the geodesic $\mathcal{A}$ by \lemref{lem:strongly-transitive}. 
This implies that every automorphism
of  $\mathcal{A}$ as a bipartite graph 
is realized by some element of $N^+$, 
showing the first claim made in part~(2)(a)(\romannumeral1). 
Since 
$T$ is the kernel 
of the action of $N^+$ 
on $\mathcal{A}$, 
this implies that 
$N^+/T$ is isomorphic to 
the group of type-preserving automorphisms 
of a triangulated line. 
Since 
the latter group 
is the infinite dihedral group, 
$D_\infty$, 
claim~(2)(a)(\romannumeral1) 
follows. 

Every element 
stabilizing 
a doubly-infinite geodesic 
and fixing an edge thereon 
must act trivially 
on that geodesic. 
This shows that 
$N^+\cap \mathcal{B}\subseteq T$. 
Since    
the opposite inclusion 
is obvious, 
we have shown that 
$N^+\cap \mathcal{B}=T$. 
In particular 
$T$ is contained in $\mathcal{B}$,
from which
the remaining claims 
of statement~(2)(a)(\romannumeral2)  
follow immediately.

Let $g$ be 
an element of $G^+$. 
By definition of $\mathcal{B}$, 
$\mathcal{B}g\mathcal{B}$ 
is the set 
of elements 
of $G^+$ 
which map 
the edge $e$ 
to some element in 
the set $\mathcal{B}g.e$. 
Since 
the group $N^+$ 
is transitive 
on the set 
of edges of $\mathcal{A}$ 
and 
the group $\mathcal{B}$ 
is transitive 
on the set of 
doubly-infinite geodesics 
containing 
the edge $e$, 
every edge 
is contained in 
the set $\mathcal{B}N^+.e$. 
We conclude that 
$G^+=\mathcal{B}N^+\mathcal{B}$. 

To finish 
the proof of part~(2)(a)(\romannumeral3),  
suppose that 
$w$ and $w'$ are 
two elements of $W$ 
such that 
$\mathcal{B}w\mathcal{B}=\mathcal{B}w'\mathcal{B}$. 
We have to show that 
$w$ equals $w'$. 
For 
every element $g$ 
of $G^+$,  
there is 
a unique edge 
of $\mathcal{A}$ 
contained in $\mathcal{B}N^+.e$. 
If $g$ is 
an element of $N^+$ 
mapping to 
an element 
$\widehat{w}$ of $W$, 
that edge is $\widehat{w}.e$. 
By the definition 
of the group $W$, 
its action 
on the set 
of edges 
of $\mathcal{A}$ 
is simply transitive, 
and 
we conclude that 
$w$ equals $w'$ 
as claimed. 
This completes the proof 
of (2)(a).

We now turn 
to part~(b) 
of claim~(2) 
which deals with 
the parametrization 
of double cosets 
in the case $G\neq G^+$. 

We have seen 
in part~(2)(a)(\romannumeral1) 
that 
all automorphisms 
of the bipartite graph $\mathcal{A}$ 
are realized by 
an element of 
the subgroup $N^+$ 
of $N$. 
Furthermore, 
of \proref{pro:rel-edgeStab}, 
guarantees 
the existence of 
an element $i$ 
in $N$ 
inverting the edge $e$. 
Therefore 
all graph automorphisms 
of $\mathcal{A}$ 
are realized by 
some element of 
the group $N$, 
showing 
the first claim 
made in part~(\romannumeral1) 
of claim~(2)(b). 
Since 
$T$ is the kernel 
of the action of $N$ 
on $\mathcal{A}$, 
this implies that 
$N/T$ is isomorphic to 
the group of automorphisms 
of a triangulated line. 
The remaining statements 
of part~(2)(a)(\romannumeral1) 
are left to the reader. 

We have seen 
in part~(\romannumeral2) of claim~(2)(a), 
which we have already established, 
that $T=N^+\cap \mathcal{B}$. 
In particular, 
the group $T$ 
is contained in $\mathcal{B}$. 
The remaining claims 
of statement~(2)(b)(\romannumeral2) 
follow 
from that fact 
and 
part~(2)(b)(\romannumeral2) 
is completely verified.

An element $g$ 
of $G$ 
which is not contained 
in $G^+$ 
can be written uniquely 
as $g_+ i$, 
where
$g_+$ is in $G^+$. 
It follows that 
$g$ is contained in 
$\mathcal{B}g_+i\mathcal{B}=\mathcal{B}g_+\mathcal{B}i\subseteq 
\mathcal{B}N_+\mathcal{B}i=\mathcal{B}N_+i\mathcal{B}\subseteq
\mathcal{B}N\mathcal{B}$,  
and 
we conclude that 
$G=\mathcal{B}N\mathcal{B}=\mathcal{B}\widetilde{W}\mathcal{B}$. 
It remains to show that 
the union $\bigcup_{\widetilde{w}\in\widetilde{W}} \mathcal{B}\widetilde{w}\mathcal{B}$ 
is disjoint. 
Now 
each element in $\mathcal{B}\widetilde{w}\mathcal{B}$ 
with $\widetilde{w}$ in $\widetilde{W}$ 
is type-preserving 
as an automorphism 
of $\mathcal{T}$ 
if and only if 
$\widetilde{w}$ is type-preserving 
as an automorphism 
of $\mathcal{A}$. 
Assume now that 
$\widetilde{w}$ and $\widetilde{w}'$ 
in $\widetilde{W}$ 
are such that 
$\mathcal{B}\widetilde{w}\mathcal{B}=\mathcal{B}\widetilde{w}'\mathcal{B}$. 
If either of 
$\widetilde{w}$ and $\widetilde{w}'$ 
is type-preserving, 
so is the other, 
both elements 
are seen to be 
contained in 
the subgroup $W$ 
of $\widetilde{W}$ 
and 
they are equal 
by part~(2)(a)(\romannumeral3), 
which we already established. 
So 
the only remaining case 
to consider 
is the one where
both $\widetilde{w}$ and $\widetilde{w}'$ 
do not preserve types. 
In this case 
choose 
elements $w$ and $w'$ 
in $W$ 
such that 
$\widetilde{w}=w.iT$
and 
$\widetilde{w}'=w'.iT$. 
Our assumption 
$\mathcal{B}\widetilde{w}\mathcal{B}=\mathcal{B}\widetilde{w}'\mathcal{B}$ 
then implies that 
$\mathcal{B}w\mathcal{B}i=\mathcal{B}w'\mathcal{B}i$, 
which is equivalent to 
$\mathcal{B}w\mathcal{B}=\mathcal{B}w'\mathcal{B}$. 
The latter condition 
implies 
$w=w'$ 
by part~(2)(a)(\romannumeral3) 
which we already established. 
We conclude that 
$\widetilde{w}$ equals $\widetilde{w}'$ 
in that case also. 
This shows that 
$\bigcup_{\widetilde{w}\in\widetilde{W}} \mathcal{B}\widetilde{w}\mathcal{B}$ is 
a disjoint union. 

In order 
to finish  
the proof 
of claim~(2) 
of Theorem~\ref{thm:rel-edgeStab} 
we need 
to establish 
the claims 
in its last paragraph. 
Only 
the statements 
about the factorization 
of the elements $\Delta_{\widetilde{w}}$ 
for $\widetilde{w}$ in $\widetilde{W}$ 
need to be proved, 
because 
the other statements 
follow 
from what 
we have 
already established 
so far.


We begin by proving 
the claims 
on $\Delta_w$  
for $w$ in $W$ 
by induction 
on the length 
of the reduced decomposition 
of $w$ 
as a word 
in the set 
of generators 
$\{s,t\}$. 
Both of these claims 
are clearly true 
if 
$w$ is the identity 
of $W$. 
Assume 
the induction hypothesis 
for elements $w'$ 
with reduced decomposition 
of length $n$ 
and 
let $w$ 
be an element 
of $W$ 
with reduced decomposition 
$s_1\cdots s_ns_{n+1}$ 
of length $n+1$. 

The word $s_1\cdots s_n$ 
is a reduced decomposition 
of $w':=s_1\cdots s_n$. 
By the induction hypothesis, 
we have 
\[
\mathcal{B}s_1\mathcal{B}s_2\mathcal{B}\cdots \mathcal{B}s_n\mathcal{B}=\mathcal{B}w'\mathcal{B}
\qquad\text{and}\qquad
\prod_{i=1}^n R(\mathcal{B}s_i\mathcal{B})= R(\mathcal{B}w'\mathcal{B})\,.
\]
We claim that 
\begin{equation}\label{eq:mult(B-doublecoset(s))}
\mathcal{B}w'\mathcal{B}s_{n+1}\mathcal{B}=\mathcal{B}w's_{n+1}\mathcal{B}
\qquad\text{and}\qquad 
R(\mathcal{B}w's_{n+1}\mathcal{B})=R(\mathcal{B}w'\mathcal{B})\,
R(\mathcal{B}s_{n+1}\mathcal{B})\,.
\end{equation}
Since $R$ is a ring homomorphism 
from the Hecke algebra 
into the integers 
the statement 
needed for 
the inductive step 
will follow 
from these equations 
and 
the induction hypothesis 
and our claims 
involving 
the elements $\Delta_w$ 
will be proved. 

Choose elements 
$n_{w'}$ and $n_{s_{n+1}}$ 
in $N^+$ 
mapping to $w'$ 
and $s_{n+1}$ respectively.
To show 
the first equation 
listed in~(\ref{eq:mult(B-doublecoset(s))}) above, 
it suffices to prove that 
$n_{w'}\mathcal{B}n_{s_{n+1}}\subseteq \mathcal{B}w's_{n+1}\mathcal{B}$, 
because 
this implies that 
$\mathcal{B}w'\mathcal{B}s_{n+1}\mathcal{B}\subseteq\mathcal{B}w's_{n+1}\mathcal{B}$; 
the opposite inclusion 
$\mathcal{B}w'\mathcal{B}s_{n+1}\mathcal{B}\supseteq\mathcal{B}w's_{n+1}\mathcal{B}$
is trivial.
In order to verify 
the above inclusion 
we apply 
a product $n_{w'}bn_{s_{n+1}}$ 
with $b$ in $\mathcal{B}$ 
to the edge $e$ 
and find that  
it sends $e$ 
to an edge 
with distance $n+1$ 
from $e$ 
on the same side of $e$ 
as the edges in $\mathcal{B}w's_{n+1}\mathcal{B}$. 
This is enough 
to imply 
the desired inclusion. 

In order 
to establish the second equation 
listed in~(\ref{eq:mult(B-doublecoset(s))}) above 
note that 
there are $R(\mathcal{B}s_{n+1}\mathcal{B})$ ways 
to extend a geodesic 
from  $e$ to $w'.e$ 
to a geodesic of length $n+1$. 

The statement involving 
the multiplication 
by the element $\Delta_i$ 
in the case $G\neq G^+$ 
follows from 
what we have just seen 
and 
$R(\Delta_i)=1$. 
Also, 
since we have 
$\Delta_i\Delta_s\Delta_i=\Delta_t$, 
we can 
omit either $\Delta_s$ or $\Delta_t$ 
from the set $\{\mathcal{B}, \Delta_s, \Delta_t, \Delta_i\}$ 
of generators 
of $\mathbb{C}[G,\mathcal{B}]$. 
The proof 
of claim~(2) 
of Theorem~\ref{thm:rel-edgeStab} 
is complete.

We now derive 
the structure 
of the Hecke algebra $\mathbb{C}[G^+,\mathcal{B}]$ 
as claimed 
in part~(3)(a) 
of Theorem~\ref{thm:rel-edgeStab}. 
We first 
verify the relations 
listed 
at the end 
of claim~(3)(a). 

For $r$ in $\{s,t\}$ 
one verifies that 
$
\mathcal{B}r\mathcal{B}r\mathcal{B}=\mathcal{B}\cup \mathcal{B}r\mathcal{B}
$.  
Hence we have 
$
\Delta_r^2=\lambda+\mu \Delta_r
$ 
with some positive integers $\lambda$, $\mu$. 
Furthermore, 
$\lambda$ is given by 
\[
\lambda=
|\mathcal{B}\backslash \mathcal{B}r^{-1}\mathcal{B}\cap \mathcal{B}r\mathcal{B}|=
|\mathcal{B}\backslash \mathcal{B}r\mathcal{B}|=:
q_r\,,
\]
and the value of $\mu$ is obtained 
by applying 
the homomorphism $R$ 
to the equality 
$
\Delta_r^2=\lambda+\mu \Delta_r
$; 
we get 
$q_r^2=\lambda+\mu q_r$. 
We deduce 
$\mu=q_r-1$, 
which proves 
the relation claimed. 

That the relations 
give a defining set of relations 
can be proved 
in the same fashion as 
in the proof of 
Theorem~3.5 in~\cite{IHES_25_5_0}. 
The rest 
of claim~(3)(a) 
is a consequence 
of what 
we have already seen. 
More precisely, 
the last two types of equations listed 
follow from 
the decomposition 
of $\Delta_w$ 
as a product 
in $\Delta_s$ and $\Delta_t$ 
corresponding to 
the reduced decomposition 
of $w$ 
established 
in part~(2)(a), 
which we already proved. 
The other types 
of equations 
follow 
using 
the following four observations: 
(\romannumeral1) under the conditions stated, 
      the reduced decomposition of $w$ 
      begins / ends with $r$; 
(\romannumeral2) using part~(2)(a),  
      we may write $\Delta_w$ as a product 
      involving $\Delta_r$ 
      at the beginning / end;  
(\romannumeral3) we may then apply the relation   
      in part~(3)(a), 
      which we 
      have already verified; 
(\romannumeral4) the desired equation 
     then follows 
     using part~(2)(a) 
     once more. 
This finishes 
the proof 
of claim~(3)(a). 

Before 
establishing 
claim~(3)(b), 
we note that 
the isomorphism 
of the algebra $\mathbb{C}[G^+,\mathcal{B}]$ 
with $\mathbb{C}W\cong\mathbb{C}D_\infty$ 
can be completed 
following the outline 
in the second paragraph 
following 
the statement 
of Theorem~1.11 in~\cite{HeckeAs+char(parab-type)(finGs:BN-pairs)}, 
thus establishing 
the part 
of claim~(1) 
concerned with $\mathbb{C}[G^+,\mathcal{B}]$. 

We turn to 
the proof of 
part~(3)(b). 
The structure of $\mathbb{C}[G,\mathcal{B}]$ 
can be described 
in terms of 
the structure of $\mathbb{C}[G^+,\mathcal{B}]$ 
using 
the reasoning 
of the analogous result, 
Proposition~3.8
of \cite{IHES_25_5_0}. 
The description 
we obtain  
will be in terms of  
the concept of 
twisted tensor product 
that is introduced in 
the following definition. 

\begin{definition}\label{def:twisted-tensor-prod}
Let $\Omega$ be 
a group 
acting by automorphisms 
on a ring $R$. 
Then the $\mathbb{Z}$-module $\mathbb{Z}[\Omega]\otimes_\mathbb{Z} R$ 
with the multiplication 
given on elementary tensors by 
\[
(\omega\otimes r)\cdot (\omega'\otimes r'):= \omega\cdot\omega'\otimes {\omega'}^{-1}(r)\cdot r'
\]
for $\omega$, $\omega'$ in $\Omega$ and $r$, $r'$ in $R$ 
is called \emph{the twisted tensor product} 
of  the given action of $\Omega$ on $R$. 
The twisted tensor product 
of an action, 
$\alpha$ say, 
of a group $\Omega$ on $R$ 
will be written 
$\mathbb{Z}[\Omega]\otimes_\alpha R$. 
\end{definition}

We return 
to the proof 
of part~(3)(b) 
of Theorem~\ref{thm:rel-edgeStab}. 

The map 
sending $i$ in $\mathbb{Z}/2\mathbb{Z}$ 
to $\Delta_i$ and the identity to $\Delta_e$ 
defines 
an injective ring homomorphism 
from $\mathbb{Z}[\mathbb{Z}/2\mathbb{Z}]$ 
into $\mathbb{Z}[G,\mathcal{B}]$ 
with image $\mathbb{Z}[N_G(\mathcal{B}),\mathcal{B}]$, 
since $\Delta_i\Delta_{\widetilde{w}}=\Delta_{i\widetilde{w}}$ 
for all $\widetilde{w}$ in $\widetilde{W}$ 
by claim~(2) 
of Theorem~\ref{thm:rel-edgeStab}, 
which we already proved. 
We will 
identify $\mathbb{Z}[\mathbb{Z}/2\mathbb{Z}]$  
with its image 
under this homomorphism. 

Part~(2)(b)(\romannumeral3) 
of Theorem~\ref{thm:rel-edgeStab} 
together with  
$\mathbb{Z}[\mathbb{Z}/2\mathbb{Z}]\cdot W=\widetilde{W}$ 
and 
$\mathbb{Z}[\mathbb{Z}/2\mathbb{Z}]\cap W=\{1\}$
identify 
$\mathbb{Z}[G,\mathcal{B}]$ 
as a $\mathbb{Z}$-module 
with the tensor product 
$\mathbb{Z}[N_G(\mathcal{B}),\mathcal{B}]\otimes_\mathbb{Z} \mathbb{Z}[G^+,\mathcal{B}]$ 
by the map $\rho\otimes \Delta_w\mapsto \Delta_\rho\Delta_w$ 
($\rho\in \mathbb{Z}[\mathbb{Z}/2\mathbb{Z}]$, $w\in W$). 

Now 
for every $\rho\in \mathbb{Z}/2\mathbb{Z}$, 
the element $\Delta_\rho$ 
is invertible 
with inverse $\Delta_{\rho^{-1}}$. 
Hence 
$\mathbb{Z}[\mathbb{Z}/2\mathbb{Z}]$ acts on $\mathbb{Z}[G^+,\mathcal{B}]$ 
through the setting 
$\rho(\Delta_w):=\Delta_\rho\Delta_w\Delta_{\rho^{-1}}=\Delta_{\rho w\rho^{-1}}$ 
($\rho\in \mathbb{Z}[\mathbb{Z}/2\mathbb{Z}]$, $w\in W$). 
Also, 
we have that 
the non-trivial element of $\mathbb{Z}[\mathbb{Z}/2\mathbb{Z}]$
sends $\Delta_w$ 
to $\Delta_{\overline{w}}$ 
for $w$ in $W$. 
Thus 
the multiplication law 
in the tensor product 
$\mathbb{Z}[\mathbb{Z}/2\mathbb{Z}]\otimes_\mathbb{Z} \mathbb{Z}[G^+,\mathcal{B}]$ 
is given by 
\[
(\rho\otimes \Delta_w)\cdot (\rho'\otimes \Delta_{w'})= 
\rho\cdot\rho'\otimes {\rho'}^{-1}(\Delta_w)\cdot\Delta_{w'}\,,
\]
for $\rho$, $\rho'$ in $\mathbb{Z}/2\mathbb{Z}$ and $w$, $w'$ in $W$.  
Hence 
$\mathbb{Z}[G,\mathcal{B}]$ exhibits 
the announced structure 
of  twisted tensor product. 
That 
the relations listed 
in part~(3)(b) 
define a presentation 
of $\mathbb{Z}[G,\mathcal{B}]$ 
follows from 
this structure result 
and 
the presentation 
of $\mathbb{Z}[G,\mathcal{B}]$ 
listed in claim~(3)(a). 
Likewise, 
the outstanding statement 
in claim~(1) 
follows. 

The proof 
of Theorem~\ref{thm:rel-edgeStab} 
is complete. 
\end{proof}

\begin{corollary}
The complete multiplication tables 
for the canonical generators 
are given by:  
\begin{enumerate}
\item 
For the algebra $\mathbb{C}[G^+,\mathcal{B}]$,  
let $w$ and $w'$ be in $W$. 
Denote by $w_{(i)}$ 
(respectively $w_{[i]}$)
the word 
in $\{s,t\}$ 
consisting of 
(respectively 
missing) 
the last 
$i$ letters 
of $w$ 
(so that 
$w=w_{[i]}w_{(i)}$) 
and 
by 
$_{(i)}w'$ 
(respectively 
${_{[i]}w'}$)
the word 
consisting of 
(respectively 
missing)
the first $i$ letters 
of $w'$ 
(so that 
$w={_{(i)}w'}{_{[i]}w'}$);  
also, 
whenever $r$ 
is in the set $\{s,t\}$, 
denote by $u$ 
the other element 
of that set; 
for a word $w$ 
with reduced decomposition $s_1\cdots s_n$  
in $W$ 
put $q_w:=\prod_{i=1}^n q_{s_i}$ 
and  
for $i$ even 
respectively 
$i$ odd 
let $q_i$ equal  
$q_r$ 
respectively 
$q_u$. 
With
this notation, 
we have
\begin{enumerate}
\item 
Suppose
the last letter 
of the reduced decomposition 
of~$w$ 
differs from 
the first letter
of the reduced decomposition 
of ~$w'$;
then 
$\Delta_w\Delta_{w'}=\Delta_{ww'}$; 
\item 
Suppose now
the last letter 
of the reduced decomposition 
of~$w$ 
coincides with
the first letter
of the reduced decomposition 
of ~$w'$ , let
$m$ be
the length 
of the shorter word 
among $w$ and $w'$ 
and let 
$s_i$ be 
the $i$-th letter 
of $w'$; 
then \\
$ 
\Delta_w\Delta_{w'}=
q_w\Delta_{w_{[m]}{_{[m]}w'}}+
\sum_{i=0}^{m-1} q_{w_{(i)}} (q_{i}-1)\Delta_{w_{[i]} s_i  {_{[i]}}w'}
$. 
\end{enumerate} 
\item 
When $G\neq G^+$, 
the products of 
elements of 
the standard basis 
for the algebra $\mathbb{C}[G,\mathcal{B}]$ 
are determined 
using 
the information 
given in~(4)(a) above 
and 
the following relations 
for $w$, $w'$ in $W$: 
$\Delta_{iw}\Delta_{w'}=\Delta_i(\Delta_w\Delta_{w'})$,
$\Delta_{w}\Delta_{iw'}=\Delta_i(\Delta_{\overline{w}}\Delta_{w'})$ 
and 
$\Delta_{iw}\Delta_{iw'}=\Delta_{\overline{w}}\Delta_{w'}$.
\end{enumerate}

\end{corollary}

\begin{proof}
Part~(1) 
is by induction 
on the quantity $m$ 
while 
part~(2) 
follows 
from 
part~(1)
and the relations 
involving 
the element $\Delta_i$ 
that were 
obtained 
in parts~(2) and~(3)(b) of \thmref{thm:rel-edgeStab}. 
\end{proof}
We close  
this section 
with an example 
showing that 
the short exact sequence 
$1\to G^+\to G\to \mathbb{Z}/2\mathbb{Z}\to 1$ 
does not always split. 

\begin{example} 
Let $K$ be a local field 
with additive valuation $\nu$ 
and 
let $G$ 
be the group $\mathrm{GL}_2(K)$ 
acting 
in the natural way 
on the Bruhat-Tits tree 
of $\mathrm{SL}_2(K)$. 
We show that 
the group $G$ 
does not contain 
an inversion 
that is 
an involution. 

First, 
by transitivity 
of $G$ on edges,  
it suffices to show that 
in $G$ 
no inversion 
of a fixed chosen edge 
is an involution. 
An element of $G$ 
that inverts 
the edge 
fixed by 
the Iwahori subgroup of $G$ 
must be 
a monomial matrix 
of the form 
\[
\left(\begin{array}{cc}0 & b \\c & 0\end{array}\right)\qquad
\text{with}\ |\nu(b)-\nu(c)|=1\,.
\]
The square of 
a matrix 
of the above form 
is $bc$ 
times the identity matrix, which,  
by the condition on the valuations 
of $b$ and $c$, 
cannot be 
the identity. 
%
\end{example}

\section{Stabilizer of an end and HNN extensions.}
\label{sec:Hecke(Stab(end),Stab(vertex))}

%
In this section 
we consider 
another 
geometrically defined class 
of Hecke pairs 
arising from 
groups acting on 
a locally finite tree, 
$\mathcal{T}$. 
The Hecke pairs 
considered 
will consist of 
the stabilizer, 
$B$ , 
of an end of $\mathcal{T}$ 
and 
its subgroup, 
$M_0$, 
stabilizing 
a distinguished vertex. 
Again 
we will assume 
that 
the group action 
is `highly transitive' 
in a sense 
to be explained shortly. 
The assumptions imposed 
will equip 
$B$ 
with 
the structure of 
an HNN-extension 
relative to 
an endomorphism $\alpha$ 
of 
$M_0$ 
such that the index $|M_0\colon \alpha(M_0)|$ 
is finite. 
(Recall that 
this HNN-extension 
can be 
defined as 
the group 
with 
set of generators 
$M_0$ 
and 
a letter $a$ 
not in $M_0$,  
subject to 
the relations 
in $M_0$ 
and 
the additional relations 
$am_0a^{-1}=\alpha(m_0)$ 
for $m_0\in M_0$.) 
We remark that 
conversely, 
every HNN-extension 
satisfying 
the latter condition 
defines 
such a tree action, 
whose 
associated HNN-extension 
is isomorphic to 
the given one 
and 
the tree action 
can be reconstructed 
from 
its associated HNN-extension 
on its `minimal subtree'. 

The following proposition 
fixes our assumptions 
and 
derives basic properties 
of the situation under consideration; 
we leave 
its proof 
to the reader. 

\begin{proposition}\label{prop:Hecke-pairs_from_HNNexts}
Let 
$\mathcal{T}$ be 
a locally finite tree, 
$\infty$ 
an end of $\mathcal{T}$,  
$B$ be a group of automorphisms 
of $\mathcal{T}$ 
that stabilize $\infty$. 
Suppose that 
some element 
of $B$, 
$s$ say,  
acts by 
a hyperbolic isometry. 
Replacing $s$ 
by another element 
if necessary, 
we may assume 
that 
$\infty$ 
is attracting  
for $s$ 
and that 
$s$ has 
the smallest 
translation length 
amongst 
elements 
of $B$.  

Denote by $M$ 
the subset 
of elements 
of $B$ 
that have a fixed point 
in $\mathcal{T}$. 
Then 
$M$ is normal in $B$ 
with infinite cyclic quotient, 
generated by 
the image of $s$. 
%
We have 
a split exact sequence 
${1\longrightarrow M\longrightarrow B\longrightarrow \mathbb{Z}\longrightarrow 1}$; 
hence 
$B$ is 
the semidirect product 
$M\rtimes \mathbb{Z}$, 
with either $1$ or $-1$ in $\mathbb{Z}$ acting 
via conjugation  
by $s$ 
on $B$. 

Let $o$ be 
a vertex 
on the axis 
of $s$ 
and 
let $M_0$ be 
the stabilizer of $o$ 
in $B$. 
%
Then 
$M_0\subseteq M$ 
and  
$M=\bigcup_{n\in\mathbb{N}}s^{n}M_0s^{-n}$. 
%
Hence 
$B$ is isomorphic to 
the HNN-extension 
of the group $M_0$ 
with respect to 
the endomorphism  
of $M_0$ 
defined by 
conjugation with $s^{-1}$, 
which 
will be 
denoted $\alpha$. 

%
The index 
$|M_0\colon \alpha(M_0)|$ 
is finite 
and 
at most 
equal to 
$|\{x\in \mathcal{T}^0\colon d(x\udodot o)=d(s^{-1}\udodot o,o)\}|-1$. 
Therefore  
$s^{-1}$ commensurates $M_0$, and hence 
so does 
$\langle s, M_0\rangle =B$. 
Thus
$(B,M_0)$ and $(M,M_0)$ 
are Hecke pairs. 
%

\end{proposition}

In the following 
we will study 
the algebra $\mathbb{C}[B,M_0]$ 
for some choices 
of $B$.  
Again 
there are 
many examples 
of groups 
whose action 
on a tree 
satisfies 
the assumptions 
of Propositon~\ref{prop:Hecke-pairs_from_HNNexts}: 
The group of 
all automorphisms 
of a locally finite (semi-)homogenous tree 
that fix an end 
satisfies 
the assumptions 
on $B$ 
above,  
as do 
the stabilizers 
of an end 
inside 
any semisimple matrix group of rank 1 over a local field 
with 
the group acting 
by automorphisms 
of its Bruhat-Tits tree.  

In contrast 
the Sections~\ref{sec:Hecke_relStab(vertex)}
and~\ref{sec:Hecke_relFix(edge)} 
however, 
the degree  
of transitivity 
imposed by 
the assumptions 
of Propositon~\ref{prop:Hecke-pairs_from_HNNexts} 
is not strong enough 
to ensure 
that all 
the resulting 
Hecke-algebras 
are isomorphic; 
the cause 
for this phenomenon 
is that if $0$ 
denotes the end of the axis of $s$ 
different from $\infty$, then 
the group $B$ does not necessarily 
act doubly-transitively on the set $B.0$. 

In the following subsection 
we show
that 
the algebra $\mathbb{C}[B,M_0]$ 
can always be described 
as the crossed product of
the algebra $\mathbb{C}[M,M_0]$ 
by
the endomorphism $\alpha$.  
We also derive 
some  general properties 
of  the algebra $\mathbb{C}[M,M_0]$. 
However, we determine 
the complete structure of 
the algebra $\mathbb{C}[M,M_0]$
only for some
 special cases only,  see  subsections \ref{subsec:Hecke-elliptic_combinatorial} and 
 \ref{subsec:Hecke-elliptic_algebraic} below.

\subsection{Reduction to the centralizer of an end relative to a vertex stabilizer}
\label{subsec:reduction->Cent(end)_Stab(vertex)}

In the context 
of Proposition~\ref{prop:Hecke-pairs_from_HNNexts} 
the group $M$ 
is said 
to be 
the `centralizer' 
of the end $\infty$  
because 
it consists of 
the elements 
of $B$ 
that fix 
some representative 
of the end $\infty$ 
pointwise 
(a property 
which 
e.g. 
the non-zero powers 
of 
the element $s$ 
do not have). 
The description 
of the algebra $\mathbb{C}[B,M_0]$ 
in terms of 
the algebra $\mathbb{C}[M,M_0]$ 
that 
we promised 
above 
exhibits 
$\mathbb{C}[B,M_0]$ 
as 
the *-algebraic semigroup crossed product 
$\mathbb{C}[M,M_0]\rtimes_\alpha \mathbb{N}$
as a consequence of the results proved 
in~\cite{lac-lar} and \cite{etal}. 
We state 
the result needed 
in the special case 
considered here 
for the convenience 
of the reader 
after clarifying 
some terms 
used in the statement.\\ 

We briefly recall first the definition of
a semigroup crossed product 
and refer 
the reader 
to~\cite[p.~422]{semiG-cross-prod+ToeplitzA(nonabGs)}
for a detailed definition. 
Suppose 
$A$ is a unital $C$*-algebra 
and 
$P$ is a semigroup 
with an action $\alpha\colon P\to \End A$ 
of $P$ 
on $A$ 
by (not necessarily unital) endomorphisms. 
The \emph{crossed product} 
of $A$  by this action of $P$ 
is an algebra, 
denoted $A\rtimes_\alpha P$, 
which is  universal and minimal  
with respect to the conditions that 
$A$  embed in $A\rtimes_\alpha P$ 
via a unital homomorphism $i_A$, 
that 
$P$ embed in the isometries of $A\rtimes_\alpha P$  
via a semigroup homomorphism $i_P$, and that 
\[
i_A(\alpha_x(a))=i_P(x)i_A(a)i_P(x)^* \qquad \text{ for } x\in P \text{ and }a\in A.
\]
A similar definition applies in the category of unital *-algebras, yielding an \emph{algebraic
semigroup crossed product}.
At the purely algebraic level, a crossed product is a twisted tensor product as defined in Definition~\ref{def:twisted-tensor-prod}.

Finally, recall that a $C$*-algebra 
is \emph{approximately finite} 
if and only if 
it can 
be written as 
an inductive limit 
of finite-dimensional 
$C$*-algebras.

\begin{theorem}
\label{thm:HNN}
Let 
$B$ be 
the HNN-extension 
of a group $M_0$ 
with respect to 
an endomorphism, 
$\alpha$,   
of $M_0$ 
such that 
$|M_0\colon \alpha(M_0)|$ is finite. 
Denote 
by $B$ 
the HNN-extension 
defined by 
the endomorphism $\alpha$ 
of $M_0$, 
name 
the automorphism 
of $B$ 
induced by $\alpha$ 
again $\alpha$ 
and 
denote 
the dilated group 
$\bigcup_{n\in\mathbb{N}} \alpha^{-n}(M_0)\subseteq B$ 
by $M$.  
Then $(B, M_0)$ and $(M,M_0)$  
are Hecke pairs 
and 
there is an action $\widehat\alpha$ of $\N$ 
by injective endomorphisms 
of the Hecke algebra $\mathbb{C}[M, M_0]$ 
given 
for $A\in M_0\backslash M/M_0$
by 
the formula 
\[
\widehat\alpha(A)  = 
{| M_0: \alpha(M_0)|}^{-1}
\sum_ {B \in M_0 \backslash \alpha^{-1}(A)/ M_0} B\,.
\]
The Hecke algebra $\mathbb{C}[B, M_0]$ 
is isomorphic 
(as a unital *-algebra) 
to 
the semigroup crossed product $\mathbb{C}[M, M_0] \rtimes_{\widehat\alpha}  \N$  
via a map 
that extends 
the canonical injection 
$\mathbb{C}[M, M_0] \hookrightarrow \mathbb{C}[B, M_0]$. 
Furthermore, 
$\mathbb{C}[M, M_0]$ is 
the union of 
a directed family 
of finite dimensional *-subalgebras, 
so 
$C^*_u[M,M_0]$ exists 
and is 
approximately finite;  
the action $\widehat\alpha$ 
extends to 
the C*-algebra level 
and 
\[
C^*_u (B,M_0)\cong  C^*_u(M,M_0) \rtimes_{\widehat\alpha} \N.
\]
\end{theorem}
\begin{proof}
That 
$(B,M_0)$ and $(M,M_0)$ 
are Hecke pairs 
follows from 
Proposition~\ref{prop:Hecke-pairs_from_HNNexts} 
and 
the statement 
at the end of 
the introductory paragraph 
to Section~\ref{sec:Hecke(Stab(end),Stab(vertex))}. 
The stated description 
of the structure 
of $\mathbb{C}[B,M_0]$ 
in terms of $\mathbb{C}[M,M_0]$ 
and $\alpha$ 
follows from 
Theorem~1.9 
of \cite{lac-lar}.

Since 
$M=\bigcup_{n\in\mathbb{N}} \alpha^{-n}(M_0)$, 
the algebra $\mathbb{C}[M,M_0]$ 
is the increasing union 
of the algebras  $\mathbb{C}[\alpha^{-n}(M_0),M_0]$; 
the latter 
are finite dimensional algebras 
because 
for each 
natural number $n$  
the index 
$|\alpha^{-n}(M_0)\colon M_0|$ 
equals $|\alpha^{-1}(M_0)\colon M_0|^n$ 
and hence
is finite. 
Hence 
$\mathbb{C}[M,M_0]$ is 
approximately finite. 
As 
each $M_0$-double-coset 
in $M$  
is contained in 
a finite dimensional \hbox{*-algebra},
it has finite spectrum. 
This
implies that 
the universal Hecke C*-algebra 
$C^*_u(M, M_0)$ exists 
and 
is approximately finite. 
The final statement 
relating 
$C^*_u (B,M_0)$ and $C^*_u(M,M_0)$ 
follows from 
Theorem~1.11 
of \cite{lac-lar}.
\end{proof}

\subsection{K-theory considerations.}
\label{subsec:compute_K-Theory}
The framework 
determined by Proposition~\ref{prop:Hecke-pairs_from_HNNexts} is 
particularly suitable to explore the K-theory 
of the algebra $\mathbb{C}[B,M_0]$.  Indeed, 
since the $C$*-algebra 
of the Hecke pair $(M, M_0)$ is approximately finite, 
it has trivial $K_1$ groups, 
and the Pimsner-Voiculescu 
six term exact sequence for crossed products
(see \cite{pim-voi}) yields

\stepcounter{theorem}\begin{equation} 
\label{pvseq}
\begin{CD}
 K_0(C^* (M,M_0))   @>{\ \ (\idx - \alpha)_*\  }>>
K_0(C^* (M,M_0)) @>{\ \ \  (i)_*\ \ \  }>> K_0(C^* (M,M_0)) \rtimes_{\alpha}{\N} ) \\
 @AAA @. @VVV  \\  
K_1(C^* (M,M_0)) \rtimes_{\alpha}{\N} )   @<{\ \ \  (i)_*\ \ \ }<<
0 
  @<{ (\idx - \alpha)_*\ }<<
 0,
\end{CD}
\end{equation}

\medskip
 Thus the $K$ groups of the Hecke C*-algebra of $(B, M_0)$ are given by  the cokernel and kernel of the map 
$
  ( \idx - \alpha)_*\  : K_0(C^* (M,M_0)) \to K_0(C^* (M,M_0)):
$
\[
K_0(C^* (B,M_0)) = \coker  ( \idx - \alpha)_* \qquad  \text{ and  } \qquad K_1(C^* (B,M_0)) = \ker  ( \idx - \alpha)_* .
 \]

Recall that  $M = \bigcup_{k\in\N} \alpha^{-k}(M_0)$ and 
let $H_n :=  \mathbb{C}[\alpha^{-k}(M_0), M_0]$;
then each $H_k$ is finite dimensional and thus a direct sum of full matrix algebras 
$H_k = \bigoplus_{i=1}^{n_k} M_{d(k,i)} $. 
The system of 
inclusions $H_k \hookrightarrow H_{k+1}$ 
is encoded by 
a  Bratteli diagram; 
see~\cite{ind-lims(fin-dimC*As)} or~\cite{dav} 
for a definition 
of the latter. 

The K-theory $K_0(H_k)$ depends only on the number  $n_i$ of full matrix summands:
\[
K_0(H_k) = \Z^{n_k}
\]
and since
\[
C^*(M,M_0) = \bigcup_{n\in \N} H_n,
\] 
one has
\[
K_0(C^*(M,M_0) ) = \lim_\to (\Z^{n_k}, B_k)
\]
where the map $\Z^{n_k} \to \Z^{n_{k+1}}$ 
is encoded by 
some integer matrix $B_k$ 
for each $k\in\N$.  
We will content ourselves with this brief general introduction here and 
tackle the explicit calculation of the $K$-theory of specific examples in a future article.

Next we will determine 
the structure of $\mathbb{C}[B,M_0]$
in the special cases in which the stabilizer of the end $\infty$
is taken first in the full 
automorphism group of a locally finite homogeneous tree, \thmref{thm:struc(Hecke(fullStab(end),Stab(vertex)},  and then in an algebraic group
acting in its Bruhat-Tits tree.
We rely on Theorem~\ref{thm:HNN} mainly for motivation and guidance, 
but we give explicit calculations in both cases.

\subsection{Stabilizer of an end relative to the stabilizer of a vertex --  full $\Aut \mathcal{T}$ case.}
\label{subsec:Hecke-elliptic_combinatorial}

In this subsection we look at 
the simplest instance of the general set up 
of Proposition~\ref{prop:Hecke-pairs_from_HNNexts}, in which $B$ consists of 
all the automorphisms of the regular tree
$\mathcal{T}$ that fix the distinguished end $\infty$. 

\begin{theorem}
\label{thm:struc(Hecke(fullStab(end),Stab(vertex)}
Let $\mathcal{T}$ be 
a locally finite homogenous tree  
of degree $q+1$, with $\infty$ 
a distinguished end and $o$ a distinguished vertex
of $\mathcal{T}$. Denote by 
$B$ the stabilizer of $\infty$, by $M$ 
the centralizer of $\infty$ 
and by $M_0$ the stabilizer 
of the vertex $o$. 
Let $s \in B$ be a hyperbolic element 
of translation length~$1$ such that 
$\infty$ is attracting for $s$ and 
the axis of $s$ contains the vertex $o$. 
Let $[s] :=  M_0sM_0$, viewed as 
an element of  $\mathbb{C}[B,M_0]$.
Then we have the following descriptions of 
$\mathbb{C}[B,M_0]$ and $\mathbb{C}[M,M_0]$.
\begin{enumerate}
\item 
For each positive integer $n$ 
define $M_n$ as the set of elements 
in $M$ that fix $s^{n}\udodot o$ 
but not $s^{n-1}\udodot o$.  
The group $M$ is 
the disjoint union of the sets $M_n$ 
where $n$ is a non-negative integer 
and this partition of $M$ is the partition 
into double cosets with respect to the subgroup $M_0$. 
Furthermore 
$M_n^{-1}=M_n$ 
for any non-negative integer $n$ 
and 
$R(M_n)=L(M_n)=(q-1)\,q^{n-1}$ 
for positive $n$.
%
The multiplication on $\mathbb{C}[B,M_0]$ and $\mathbb{C}[M,M_0]$ is determined by
\begin{eqnarray}
{[s]^*}^n [s]^n &=& q^n\qquad \qquad  \text{ for } n\geq 0 \label{5.4.3.a.1}\,, \\
 {[s]}^n{[s]^*}^{n} &=& \sum_{i=0}^n M_i \qquad \qquad \text{ for } n\geq 0 \label{5.4.3.a.2}\,,\\
 M_mM_n &=& \ M_nM_m=(q-1)q^{n-1}M_{m} \qquad \qquad \qquad\text{ for } m>n>0\label{5.4.3.b.1}\,,\\
M_m^2 &=& \ (q-2)q^{m-1}M_m+(q-1)q^{m-1}\sum_{i=0}^{m-1}M_i \qquad \text{ for } m>0 \,.\label{5.4.3.b.2}
\end{eqnarray}

\item 
\label{cor:struc(C[M,M_0],C[B,M_0])}
Let $\mu:=q^{-1/2}[s]$. Then
$\mu$ is a nonunitary isometry. The algebra $\mathbb{C}[M,M_0]$ is generated by 
the range projections $\{ \mu^n {\mu^*}^n:  n\in \N\}$ and,
as such, is isomorphic to 
the algebra of eventually constant sequences 
of complex numbers; in particular, it is abelian. 
The algebra $\mathbb{C}[B,M_0]$ is generated as a *-algebra
by $\mu$ and, as such, it is isomorphic to the universal *-algebra
generated by an isometry. 
\end{enumerate}
\end{theorem}
\begin{proof}
We begin 
by proving 
part~(1), 
which identifies 
the elements of 
the canonical basis 
of the algebra $\mathbb{C}[M,M_0]$ 
in terms of 
their action 
on the tree. 

That $M$ is 
the disjoint union $\bigsqcup_{n\in\mathbb{N}} M_n$ 
follows from our initial definition of $M_0$ and 
the definition 
of the sets $M_n$ for $n>0$.
It is also clear from the  definition that 
$M_n^{-1}=M_n$ for every $n>1$
and this,  in turn, will 
imply that 
$R(M_n)=L(M_n)$ 
for any $n$,  
once we know that 
each $M_n$ 
is an $M_0$ double coset. 

Since $M_0$ fixes 
all the vertices $s^{m}\udodot o$ 
for $m$ a non-negative integer, 
it follows that $M_0M_nM_0\subseteq M_n$ 
for all non-negative integers $n$, 
and hence 
each $M_n$ 
is a union 
of $M_0$ double cosets 
in $M$. 
To see that 
each $M_n$ 
consists of 
a single double coset, 
it suffices to consider 
positive $n$. 
Suppose then that 
$m$ and $m'$ 
are both contained 
in $M_n$ 
for $n$ positive. 
Let $m_0\in M_0$ be 
an element that maps 
$m\udodot o$ into  $m'\udodot o$ 
(there is such an element $m_0$, because the path 
from $s^{n}\udodot o$ to  $m\udodot o$ is of the same length as 
that from  $s^{n}\udodot o$ to  $m'\udodot o$, 
and both paths are disjoint from  the path 
 $s^{n}\udodot o$ to  $o$).
Then the element 
$m_0':={m'}^{-1} m_0 m$ fixes $o$ 
and thus is in $M_0$. 
We conclude that $m=m_0^{-1}m'm_0'\in M_0m'M_0$ 
and thus $m$ and $m'$ are
in the same 
$M_0$ double coset, 
as claimed. 

It remains 
to verify that 
$L(M_n)=(q-1)\,q^{n-1}$ 
for positive $n$.
An $M_0$ left coset inside $M_n$ is uniquely determined by 
the image of $o$ under any of its elements. 
To count the possible images 
of $o$ under elements of $M_n$, 
observe that any vertex at distance $n$ 
from the vertex $s^{n}\udodot o$ is a possibility, 
except those in the branch starting at $s^{n-1}\udodot o$. 
Thus, there are $(q-1)$ possibilities 
for the image of $s^{n-1}\udodot o$, 
and then $q$ possibilities 
for each of the following vertices 
$s^{n-2}\udodot o$, $s^{n-3}\udodot o$, up to $o$.
Thus there are $(q-1) q^{n-1}$ 
possible images for $o$ 
under elements of $M_n$. 
The proof of part~(1) is now complete. 
As an intermediate step 
towards the multiplication formulae, we will prove 
\begin{align}
\label{one}
M_0\subseteq s^nM_0s^{-n}=&\  M_{s^n\udodot o} \qquad \text{ for } n\geq 0 \text{, and}\\
\label{two}
(M_0sM_0)^n (M_0s^{-1}M_0)^n=&\ s^nM_0s^{-n} \qquad \text{ for } n>1.
\end{align}
We begin 
by proving~(\ref{one}). 
By definition,  $M_0$ is the stabilizer $M_o$ of
of $o$ in $M$, so
$s^nM_0s^{-n}=s^nM_os^{-n}=M_{s^n\udodot o}$. 
Since $o$ lies on the axis 
of $s$ and since the fixed end $\infty$ is attracting for $s$,
every element of $M$ that fixes $o$,  
also fixes $s^n\udodot o$ for $n\geq 0$.
This implies that
$M_0=M_o\subseteq M_{s^n\udodot o}=s^nM_0s^{-n}$, finishing the proof of
(\ref{one}).
Claim~(\ref{two}) 
can now be derived 
by induction on $n$
using claim~(\ref{one}).
Since $R(s) = |M_0 / s\inv M_0 s| = q$ and since $L(s) = 1$ by 
\eqref{one} above with $n=1$, we may apply
\cite[Theorem 1.4]{bre} and conclude 
that $\mu:=q^{-1/2}[s]$ is an isometry, from which we see that
$[s]^*[s] = q 1$.
An easy induction argument then proves \eqref{5.4.3.a.2}.
Using~(\ref{two}) 
and~(\ref{one}) we see that 
the element $[s]^n[s]^{*n}$ 
of the algebra $\mathbb{C}[B,M_0]$ 
is supported on 
the set $M_{s^n\udodot o}$.
which decomposes as 
the disjoint union $\bigsqcup_{k=0}^n M_k$ 
by the definition 
of the $M_k$'s. 
Since each $M_k$ is an $M_0$ double coset,
$[s]^n[s]^{*n}$ can be expressed as 
a linear combination 
$\sum_{i=0}^n k_i(n) M_i$,
and from claim~(2)  we see that $k_i(n)>0$
for all non-negative integers $n$ and all $0\le i\le n$ .

This observation, when combined  with the homomorphism $R$ of $\Z[M,M_0]$ to $\Z$
suffices to compute the coefficients $k_i(n)$.
%
Apply 
the algebra-homomorphism $R$ 
to the equation 
\begin{equation}
\label{eqn:ss*}
[s]^n{[s]^*}^{n}=\sum_{i=0}^n k_i(n) M_i\,.
\end{equation}
We have already proved that 
 $R(M_i) = (q-1) q^{i-1}$ that 
 $R(s)=q$ and that $R(s^{-1})=L(s)=1$. 
Using this information 
when applying $R$ 
to equation~(\ref{eqn:ss*})
we get 
\begin{equation}
\label{eqn:coeff(ss*)}
q^n=k_0(n)+(q-1)\sum_{i=1}^nk_i(n)q^{i-1}\,.
\end{equation}
If 
we calculate 
the sum 
on the right hand side 
of equation~(\ref{eqn:coeff(ss*)}) 
substituting 
the minimal possible value $1$ 
for the unknown quantities $k_i(n)$ 
we obtain $q^n$ 
as the result. 
Since 
all summands 
on the right hand side 
of equation~(\ref{eqn:coeff(ss*)}) 
are positive,   
none of the numbers $k_i(n)$ 
can be bigger than $1$ 
and 
we derive that 
$k_i(n)=1$ for all non-negative integers $n$ and all $0\le i\le n$. 
This finishes the proof 
of \eqref{5.4.3.a.2}.

Since $M_n\inv = M_n$ for every $n\geq 0$,  the double cosets are self-adjoint, and
$\mathbb{C}[M,M_0]$ is abelian.
The multiplication table in equations 
\eqref{5.4.3.b.1} and \eqref{5.4.3.b.2} 
can now be obtained via a recursion using
\eqref{5.4.3.a.1} and \eqref{5.4.3.a.2}.

Applying basic linear algebra 
to the expressions $[s]^n[s]^{*n} = \sum_{i=0}^n M_i$ 
and recalling that  
of $\mu = q^{1/2} [s]$,  
we see that  each $M_k$ 
is a linear combination 
of range projections 
of powers of $\mu$. 
Since 
the $M_k$ 
for $k$ in $\mathbb{N}$ 
form a basis 
of $\mathbb{C}[M,M_0]$, 
the first statement 
of part~(1) 
follows. 
The isomorphism 
with the  algebra $C_c$ 
of eventually constant sequences 
is achieved as follows. 
Let $P_i:=\mu^i\mu^{*i}$ for each $i\geq 0$ and note that 
 that the elements $Q_i := P_i-P_{i+1}$ are self adjoint and satisfy $Q_iQ_j = \delta_{i,j} Q_i$,
 and $Q_i Pj = 0$ if $i<j$ and $1$ if $i\geq j$. Since the generators of $C_c$,
 $q_i = \delta_i$ and $p_j = \chi_{\{j, j+1, \ldots\}}$
have the same multiplication table, we 
map the sequence 
$(x_0,x_1,\ldots,x_N,x_N,\ldots )\in C_c$ 
to the combination
\begin{eqnarray*}
x_0Q_0+x_1Q_1+&\cdots& + x_N P_N \\
&=&
x_0 +(x_1-x_0)P_1+(x_2-x_1)P_2+\cdots+(x_N-x_{N-1})P_N.
\end{eqnarray*}
In order to verify that the map is isometric, note that
\begin{eqnarray*}
\| x_0 +(x_1-x_0)P_1+(x_2-x_1)P_2+&\cdots&+(x_N-x_{N-1})P_N\|  \\
&=&
\max
\left\{|x_0|,\left| x_0 +\sum_{i=1}^N x_i-x_i-1\right| 
\colon 1\leq n\leq N\right\} \\
&=& 
\max \{ |x_n|\colon 0\leq n\leq N\} \\
&=&
\| (x_0,x_1,\ldots,x_N,x_N,\ldots ) \|_{C_c}.
\end{eqnarray*}
Since 
the algebra $\mathbb{C}[B,M_0]$ 
is generated 
by $\mu$ and $\mathbb{C}[M,M_0]$,  
the second statement in part~(1)
follows from the first 
and  the universal property 
of the Toeplitz-algebra (see \cite{C*-Alg_gen(isom)}).   
This finishes the proof 
of part~(2) 
and of the theorem. 
\end{proof}

\subsection{The centralizer of an end relative to a vertex stabilizer --- algebraic group cases}
\label{subsec:Hecke-elliptic_algebraic}

In this subsection 
we derive 
some information 
on the Hecke algebra 
of the 
centralizer of 
an end 
of the Bruhat-Tits tree 
relative to a vertex stabilizer 
for 
semisimple matrix groups 
of rank~$1$ 
over a local field, 
using 
the structure 
of these groups. 
In these groups 
the subgroup 
$M$ 
decomposes as 
a semidirect product 
$Z_0\ltimes U$, 
and 
we have $M_0U=M$. 
This 
can be used 
to identify 
the Hecke algebra 
$\mathbb{C}[M,M_0]$ 
with 
a subalgebra 
of $\mathbb{C}[U,U\cap M_0]$ 
and 
to explain 
the multiplication 
of $\mathbb{C}[M,M_0]$ 
in terms of 
the action 
of $Z_0$ on $U$, 
see Theorem~\ref{thm:C[M,M_0]-via-C[U,U_0]}. 
As this subsection 
is rather technical,  
those unfamiliar with 
the underlying general theory 
might find it helpful to keep in mind 
the special case 
$G = \mathrm{SL}_2(k)$ 
discussed below in Example ~\ref{illustration}. A similar result for  matrix groups over 
algebraic number fields can be found in~\cite[Proposition 1.6 (4)]{lvf} where the structures denoted $\mathfrak A$ and $A_\theta$  are Hecke algebras, as explained
in~\cite[page 38]{lvf} which act on a countably infinite product of trees. 

\begin{theorem}\label{thm:C[M,M_0]-via-C[U,U_0]}
Let $G$ be 
a semisimple matrix group   
over a local field $k$ 
and 
$\infty$ an end 
of the Bruhat-Tits tree $\mathcal{T}$ 
of $G$ 
over $k$.  
The stabilizer $B$ 
of $\infty$ 
is 
the group 
of $k$-rational points 
of a $k$-parabolic subgroup, 
$\mathbf{B}$,  
of $G$; 
let $U$ be 
the group 
of $k$-rational points 
of 
the unipotent radical 
of $\mathbf{B}$ 
and 
$Z$ 
the group 
of $k$-rational points 
of a maximal $k$-split torus 
contained 
in $\mathbf{B}$. 
Denote by $M$  
the centralizer 
of $\infty$,  
by $M_0$ 
the stabilizer in $M$ (equivalently, in $B$)
of a distinguished vertex, 
$o$, 
of $\mathcal{T}$. 
Put $U_0:=M_0\cap U$. 
Denote 
by $0$ 
the end 
different from $\infty$ 
fixed by $Z$. 
We have 
$U\subseteq M$,  
$B=Z\ltimes U$
and 
$M=(Z\cap M)\ltimes U$. 
The group $Z\cap M$ 
fixes the line 
connecting $\infty$ to $0$ 
and 
we have $Z\cap M\subseteq M_0$,  
and 
$M_0=(Z\cap M)\ltimes U_0$, 
in particular 
$M_0U=M$. 
Also, 
the group $Z_0:=Z\cap M$ 
coincides with 
the elements of $Z$ that 
fix the point $o$. 

The following statements hold:
\begin{enumerate}
\item 
For $u$ in $U$ 
we have 
(writing conjugation 
as an exponent 
on the left):
\begin{align*}
M_0uM_0\cap U &=M_0uM_0U_0\cap U={^{M_0}u}\,U_0\\
{^{m_0}}uU_0 &={^{m'_0}}uU_0
\quad\text{if and only if}\quad
m_0{{m_0}'}^{-1}\ \text{fixes}\ u\udodot o
\end{align*}
\item 
Every element 
of the standard basis 
of $\mathbb{C}[M,M_0]$ 
can be 
represented as 
$M_0uM_0$ 
with $u$ in $U$ 
and 
for $u$, $u'$ 
in $U$ 
we have 
$M_0uM_0=M_0u'M_0$ 
if and only if 
$u'
\in {^{M_0}u}\,U_0$;  
if 
$U$ is abelian,  
we can replace 
the last condition 
by 
$u'
\in {^{Z_0}u}\,U_0$. 
\item 
The map $\nu\colon \mathbb{C}[M,M_0]\to \mathbb{C}[U,U_0]$ 
defined 
for $m$ in $M$ 
by 
\[
\nu(M_0mM_0):=
\sum_{h\colon U_0\backslash (M_0mM_0\cap U)/U_0} 
U_0hU_0 
\]
induces  
an isomorphism 
of *-algebras 
between 
$ \mathbb{C}[M,M_0]$ 
and its image 
under $\nu$. 
\item 
If 
$U$ is abelian, 
and $u$ 
is in $U$, 
we have 
\[
\nu(M_0uM_0)=
\sum_{h\colon  {^{Z_0}u}U_0/U_0} 
hU_0 
\in \mathbb{C}[U/U_0]\,.
\]
\end{enumerate}
\end{theorem}
\begin{proof}
The structural properties 
of end stabilizers 
in semisimple matrix groups 
of rank~$1$ 
over 
a local field, 
that 
are stated 
in the introductory paragraph 
can be found 
for example 
in Subsection~1.4 
in~\cite{cusps(lat<rk1.localF)}. 

The first 
displayed formula 
in claim~(1) 
is proved 
using that 
$U$ is normal in $M$. 
This fact 
together with 
the definition 
of $M_0$ 
yields 
the second displayed formula 
of claim~(1) 
and establishes 
the first claim. 

In claim~(2), 
that 
every element 
of the standard basis 
of $\mathbb{C}[M,M_0]$ 
can be 
represented as 
$M_0uM_0$ 
with $u$ in $U$ 
follows from 
$M_0U=M$, 
which is 
stated 
in the introductory paragraph 
of Theorem~\ref{thm:C[M,M_0]-via-C[U,U_0]}. 
The rest 
of claim~(2) 
follows from 
claim~(1) 
which we already proved,  
while the statement 
about the case 
where 
$U$ is abelian 
follows from 
the equation 
$M_0=Z_0U_0$
that 
is obtained 
from 
the introductory paragraph  
of Theorem~\ref{thm:C[M,M_0]-via-C[U,U_0]}. 

Claim~(3) 
is essentially 
Proposition~6.4 in~\cite{kri} 
(we only need 
that 
$(M,M_0)$ is 
a Hecke pair,  
$U$ is a subgroup 
of $M$ 
such that 
$M=M_0U$ 
and that 
$U_0=U\cap M_0$). 
Concretely, 
setting 
$G:=M$, $H:=U$ and $S:=M_0$ 
in \textit{loc.cit.\/} 
yields
the corresponding 
claim 
for the Hecke-algebras 
over the integers. 
Tensoring 
with $\mathbb{C}$ 
we obtain 
the claim 
on the level 
of algebras
over the field $\mathbb{C}$,  
and 
a direct calculation 
shows that 
$\nu$ commutes with 
the *-operation, 
thus completing
the proof 
of~(3). 

Finally, 
claim~(4) 
follows from 
follows from 
claims~(1) to~(3) 
and 
the proof 
of Theorem~\ref{thm:C[M,M_0]-via-C[U,U_0]} 
is complete. 
\end{proof}

To give 
the reader 
some impression 
on the complexity 
saved 
by applying 
the map $\nu$, 
we remark that 
the group $U$ 
is 
either abelian 
or metabelian\footnote{ 
The metabelian case 
only arises 
if 
the affine root system 
of $G$ 
is of type $\widetilde{BC}_1$.}   
and 
that 
$Z$ 
acts via characters 
on the root subgroups 
whose product 
is $U$. 
As a further, 
more concrete, 
illustration, 
we now discuss 
the case 
of the group $\mathrm{SL}_2(k)$, 
where 
$k$ is 
a local field 
with ring of integers $\mathcal{O}$; 
the reader 
may prefer 
to specialize 
further 
and 
assume that 
$k=\mathbb{Q}_p$ 
and 
$\mathcal{O}=\mathbb{Z}_p$.

\begin{example}[Illustration of $(B,M,M_0,U,U_0,Z,Z_0)$ \label{illustration}
in the case of $\mathrm{SL}_2(k)$, 
$k$ a local field]
\label{ex:(B,N,U,Z,Z_0)[SL_2]}
~\\
The groups $B$, $M$, $M_0$, $U$, $U_0$, $Z$ and $Z_0$ 
depend on 
the choice of 
two ends 
and 
a vertex 
on the line 
connecting them 
in the tree $\mathcal{T}$.  
The ends 
and the vertex 
defined by 
the standard basis 
in the model of 
the Bruhat-Tits-tree $\mathcal{T}$ 
of $\mathrm{SL}_2(k)$ 
as outlined in~\cite{ser-tre} 
leads to
the following 
list: 
 \begin{enumerate}
\item 
$B$ is 
the group 
of upper-triangular matrices 
in $\mathrm{SL}_2(k)$; 
\item 
$M$ is 
the subgroup 
of $B$,  
whose 
diagonal entries 
are contained in $\mathcal{O}$; 
\item 
$M_0$ is 
the group 
of upper-triangular matrices 
in $\mathrm{SL}_2(\mathcal{O})$; 
\item 
$U$ is 
the group 
of  upper-unitriangular matrices 
in $\mathrm{SL}_2(k)$ (so all diagonal entries 1); 
\item 
$U_0$ is 
the group 
of upper-unitriangular matrices 
in $\mathrm{SL}_2(\mathcal{O})$ (so all diagonal entries 1); 
\item 
$Z$ is 
the group 
of diagonal matrices 
in $\mathrm{SL}_2(k)$; 
\item 
$Z_0$ is 
the group 
of diagonal matrices 
in $\mathrm{SL}_2(\mathcal{O})$. 
\end{enumerate}
The group $U$ 
is isomorphic to 
the additive group 
of $k$
and 
$U_0$ is isomorphic to 
the additive group 
of $\mathcal{O}$; 
hence 
$\mathbb{C}[U,U_0]\cong\mathbb{C}[k/\mathcal{O}]$. 
Also, 
$Z$ is isomorphic to 
the multiplicative group 
of non-zero elements in $k$, 
while 
$Z_0$ is isomorphic to 
the group of units 
of $\mathcal{O}$; 
conjugation by 
an element $\lambda\in k^\times\cong Z$  
maps 
an element $u\in k\cong U$ 
to $\lambda^2u$. 
Hence, 
in this case, 
we have: 
\[
\nu\left(M_0\left[\begin{array}{cc}1 & u \\0 & 0\end{array}\right]M_0\right)=
\sum_{h\colon  {(\mathcal{O}^\times)^2u}+\mathcal{O}/\mathcal{O}} 
h\mathcal{O} 
\in \mathbb{C}[k/\mathcal{O}]\,.
\]
\end{example}


\end{document}